\documentclass[11pt,a4paper]{article}
\pdfoutput=1
\usepackage{amssymb}  
\usepackage{amsthm}
\usepackage{amsmath}
\usepackage{amsfonts}
\usepackage{moreverb}
\usepackage{graphicx}
\usepackage{float}
\usepackage{enumitem}
\usepackage{caption}
\usepackage{subcaption}
\usepackage{authblk}

\title{Active Fault Tolerant Flight Control System Design}
\author[1]{Rudaba Khan}
\author[2]{Paul Williams}
\author[2]{Paul Riseborough}
\author[1]{Asha Rao}
\author[1]{Robin Hill}
\affil[1]{Department of Mathematics and Geospatial Science, RMIT University, Melbourne, Australia.}
\affil[2]{BAE SYSTEMS Australia, Melbourne Australia.}

\begin{document}
\maketitle
\begin{abstract}
In this paper we investigate the design of an active fault tolerant control system applicable to autonomous flight.  The system comprises a nonlinear model predictive based controller integrated with an unscented Kalman filter for fault detection and identification.  We apply the fault tolerant control system design to a generic aircraft model, and simulate a failed engine scenario.  The results show that the system correctly identifies the fault within seconds of occurrence and updates the nonlinear model predictive controller which is then able to reallocate control authority to the healthy actuators based upon up to date fault information.  
\end{abstract}

\maketitle

\section{Introduction}
Fault tolerant flight control is a very active area of research, however the work has concentrated on large manned aircraft with little work being done for unmanned air vehicles. In this paper, we fill this gap by developing an active fault tolerant control system for autonomous flight consisting of a nonlinear model predictive control (NMPC) based controller integrated with a fault detection and identification (FDI) system. The FDI developed here uses an unscented Kalman filter (UKF) that actively seeks to predict faults in the system and provides parameter updates to the NMPC controller for reconfiguration.  We successfully demonstrate the ability of our design to identify engine failure and use the healthy actuators to continue the mission. 

Most of the FTC schemes currently in use consist of a combination of two or more control methods.  For example the FTC system developed by the Intelligent Flight Control System (IFCS) F-15 program at NASA (\cite{perhinschi2004performance}, \cite{perhinschi2006comparison}) is based on nonlinear dynamic inversion augmented with a neural network to compensate for inversion errors and changes in aircraft dynamics due to damage or failure of a primary control surface.  Another example of combined methods for FTC is work conducted by Shin and Gregory \cite{shin2007robust}, with the FTC method based on robust gain scheduling (GS) control concepts using a linear parameters varying (LPV) control synthesis for civil transport aircraft.  Yang and Lum's \cite{yang2003fault} solution, tested on simulation models of the F-16 aircraft with stuck actuator faults, bases the FTC on $H_\infty$ and peak-to-peak gain performance indices in a multiobjective optimisation setting with the algorithms based on linear matrix inequalities (LMIs).  The active FTC method developed by Yu and Jiang \cite{Yu2012Impairments} models the impairments as a polytopic linear paramater varying (LPV) system.  Ye et. al \cite{Ye2010Uncertain} use a linearised model of the F-18 (longitudinal motion only) with the FTC based on $H_\infty$ in an LMI framework similar to Yang and Lum \cite{yang2003fault}.  

Our FTC system is based on only one method, nonlinear MPC, which is used in both the fault free and the fault cases.  Due to the inherent fault tolerant capabilities of NMPC \cite{maciejowski1998implicit}, the system is capable of control configuration without the need for another control scheme.  The model developed for this work is detailed in section \ref{section:modelDescription}.  This is followed by a description of the NMPC based FTC system in section \ref{section:FTC}.  In this section the simulations assume that FDI information is provided and the results show that our design is able to reallocate control distribution in the event of an actuator failure.  Section \ref{section:FDI} details the development of a UKF based FDI scheme.  The results of this section reveal that there are many issues relating to observability which need to be addressed.  This section identifies areas for further research.  The order of the model is reduced in section \ref{section:engineFailure} and the fault simulated is loss of engine failure.  The results show that our FTC system design successfully identifies the fault and redistributes control allocation to the healthy actuators.

\section{Model Description}\label{section:modelDescription}
\begin{figure}[H]
\center
\includegraphics[scale=0.3]{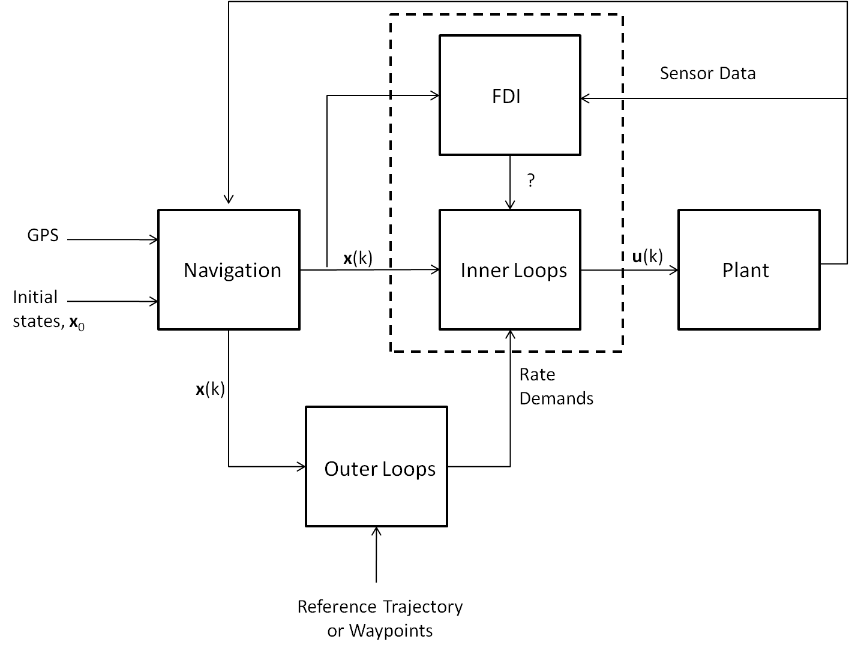}
\caption{Aircraft System with FTC} 
\label{fig:chap5_blockDiag_withFDI}
\end{figure}

The system given in figure \ref{fig:chap5_blockDiag_withFDI} was modelled in MATLAB/Simulink with a full 6DoF aircraft model.  

The plant model used here doubles as the navigation model with the guidance system provided by Williams \cite{PWOuterLoopModel}.  The Dryden Wind Turbulence model from the Aerospace Blockset \cite{MATLAB_Aero} was used to model wind and turbulence.  The guidance subsystem is supplied with a series of way points and provides the controller with angular rate information calculating the angular rates required to maintain the reference path.  The inner loops consist of two controllers, an NMPC controller to control angular rates and a speed control loop, a simple PI controller, to maintain a desired speed.  An integrator is also implemented to calculate the integrated errors in the angular rates.

\begin{figure}[H]
%\hspace*{-0.5in}
\center
\includegraphics[scale = 0.3]{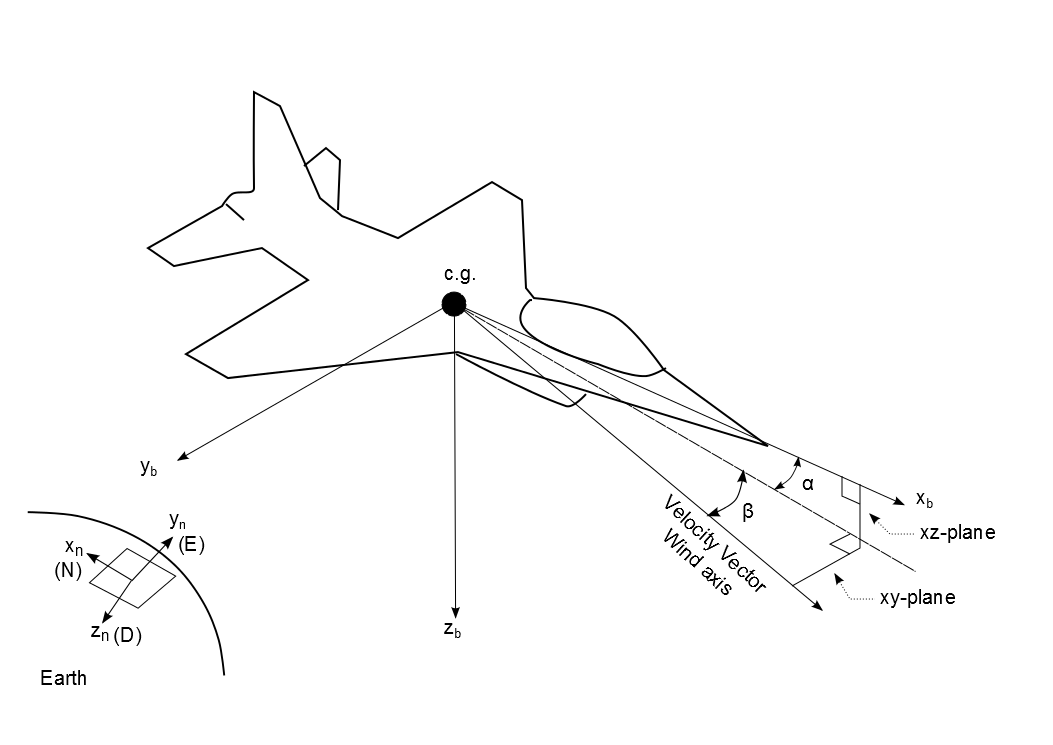}
%\vspace{-0.5in}
\caption{Aircraft Coordinate Frames}
\label{fig:chap5_axes} 
\end{figure}

The state vector of the aircraft plant model is:

\begin{equation}\label{eqn:chap5_plantStateVec}
\mathbf{x} = \left[x_N \;\; x_E \;\; x_D \;\; V_N \;\; V_E \;\; V_D \;\; \phi \;\; \theta \;\; \psi \;\; p \;\; q \;\; r\right]^\intercal,
\end{equation}

where $x_N, x_E, x_D$ are north, east and down position coordinates respectively, given in an earth fixed tangent frame, called the navigation frame (or NED frame), denoted by the subscript \textbf{n}, (see figure \ref{fig:chap5_axes}).  The NED is an Earth fixed frame, with the origin located at a point on the Earth.  In practice this origin is defined at the point where the aircraft is initialised for flight.  The vectors $V_N, V_E, V_D$ are the components of the velocity vector in the north, east, and down directions respectively, $\phi, \theta, \psi$ are the aircraft orientation angles roll, pitch and yaw respectively and, $p,q,r$ are the roll, pitch and yaw angular rates respectively.

A specific aircraft model is required for the development of a fault tolerant controller.  The generic aircraft model developed here for control law design and validation is based on the McDonnell Douglas F-4 aircraft \cite{garza2003collection}.   

The generic aircraft model developed for this work has the same aerodynamic characteristics as the F-4 aircraft but with the following mass and size properties: Wing Area $S=20m$, Mean Aerodynamic Chord $\bar{c}=3m$, C.G location $x_{c.g}=0$, C.G reference location $x_{c.g.ref}=0$, Weight 1,177 (kg), $I_X = 2,257 \text{kg.m}^2$), $I_Y =11,044 \text{kg.m}^2$, $I_Z= 12,636 \text{kg.m}^2$) and $I_{XZ}=106\text{kg.m}^2$. 

The expressions for the non-dimensional force and moment equations are given in the appendix along with the thrust model used.  The controls for the aircraft model are throttle ($\delta_{th}$), elevator ($\delta_e$), aileron ($\delta_a$) and rudder ($\delta_r$), with the control vector given by:

\begin{equation}
\mathbf{u} = \left[\delta_{th},\,\,\,\delta_a,\,\,\,\delta_e,\,\,\,\delta_r\right]^\intercal.
\end{equation}

\section{Fault Tolerant Flight Control}\label{section:FTC}
The state vector of the NMPC controller is given by:

\begin{equation}\label{eqn:xnmpc}
\mathbf{x}_{\text{nmpc}}=\left[p,\,\,q,\,\,r,\,\,I_{p},\,\,I_{q},\,\,I_{r},\,\,\delta_e,\,\,\delta_a,\,\,\delta_r,\,\,\Delta\delta_e,\,\,\Delta\delta_a,\,\,\Delta\delta_r\right]^\intercal,
\end{equation}

where $p$, $q$ and $r$ are the roll rate, pitch rate and yaw rate respectively, $I_{p}$, $I_{q}$ and $I_{r}$ are respectively, the integrated errors in $p$, $q$ and $r$ used to minimise the steady state errors, $\delta_e$, $\delta_a$ and $\delta_r$ are the elevator, aileron and rudder deflections respectively and $\Delta\delta_e$, $\Delta\delta_a$ and $\Delta\delta_r$ are the elevator, aileron and rudder rates.

The prediction model of the nonlinear MPC controller is as follows:
\begin{eqnarray}
\dot{p} &=& \left(c_1\,r+c_2\,p+c_4\,h_{eng}\right)q + \bar{q}\,S\,b\left(c_3\,C_l + c_4\,C_n\right),\\
\dot{q} &=& \left(c_5\,p - c_7\,h_{eng}\right)r-c_6\left(p^2-r^2\right) + \bar{q}\,S\,\bar{c}\,c_7\,C_m,\\
\dot{r} &=& \left(c_8\,p-c_2\,r+c_9\,h_{eng}\right)q + \bar{q}\,S\,b\left(c_4\,C_l+c_9\,C_n\right),\\
\dot{I}_p &=& \hat{p} - p_{dem},\\
\dot{I}_q &=& \hat{q} - q_{dem},\\
\dot{I}_r &=& \hat{r} - r_{dem},
\end{eqnarray} 

where $\hat{p}$, $\hat{q}$ and $\hat{r}$ are the predicted angular rates and $p_{dem}$, $q_{dem}$ and $r_{dem}$ are the demanded angular rates.  The terms $c_1$ to $c_7$ are the moments of inertia as defined in \cite{garza2003collection}.  $h_{eng}$ is the distance of the engine from the center of gravity and is taken to be 0m and $\bar{q}$ is the dynamic pressure.  $C_l$, $C_m$ and $C_n$ are the non-dimensional moment coefficients and $C_X\,C_Y$ and $C_Z$ are the non-dimensional force coefficients.

A pseudospectral discretisation method is used \cite{RKPaper1} and the following control problem is solved at each time step:

\begin{equation}
\begin{split}
\min_{\mathbf{x},\mathbf{u}}\, \frac{H_p}{2}\;\sum_{j = 1}^{j = N+1} \bigg(&\big\Vert \mathbf{\omega}(j) - \mathbf{\omega}_{\text{dem}}(j)\big\Vert_{Q_\omega}^2 +\big\Vert \bold{{\dot{I}}}(j)\big\Vert_{Q_I}^2\\
& + \big\Vert \Delta\mathbf{u}(j)\big\Vert_{Q_{u}}^2 + \big\Vert \mathbf{a}_N(j)\big\Vert_{Q_{a}}^2\bigg)\;w(j),
\end{split}
\end{equation}

subject to
\begin{eqnarray}
\left(\frac{t_f-t_0}{2}\right)\mathbf{D}_{j,k}\mathbf{x}_j - \mathbf{\dot{x}}_j &=& 0, \\
\mathbf{\omega}(j_0) - \mathbf{\omega}_{\text{dem}}(j_0) &=& 0,\\
\mathbf{x}_{lb}  \leq    \mathbf{x}  \leq  \: \mathbf{x}_{ub},\\
\mathbf{u}_{lb}  \leq    \mathbf{u}  \leq  \: \mathbf{u}_{ub},\\
\Delta\mathbf{u}_{lb}  \leq    \Delta\mathbf{u}  \leq  \: \Delta\mathbf{u}_{ub}, \label{eq:chap5_6DOF_cons}
\end{eqnarray}
where $\mathcal{D}_N$ is a spectral differentiation matrix \cite{RKPaper1}, $N$ refers to the number of discretisation (or coincidence) points, $t_0$ and $t_f$ are the initial and final times of the prediction horizon window and the term $\bold{I}$ is the vector of integrated errors .  The state vector $\mathbf{x}$ is defined in \eqref{eqn:xnmpc}.  $\Delta\mathbf{u}$ are the control input rates and $\mathbf{a}_N$ are the accelerations in the navigation frame.  The constraints on the states: $p$: None, $q$: None, $r$: None, $I_p$: None, $I_q$: None, $I_r$: None, $d_A$: $\pm 20\deg/\sec$, $d_E$: $\pm 20\deg/\sec$, $d_R$: $\pm 20\deg/\sec$, $\Delta d_A$: $\pm 200\deg/\sec$, $\Delta d_E$: $\pm 200\deg/\sec$ and $\Delta d_R$: $\pm 200\deg/\sec$.  $Q_\omega$, $Q_I$, $Q_{u}$ and $Q_{a}$ are diagonal weighting matrices with the following values along the diagonals 1, 0.000001, 0.05 and 1 respectively.  These weighting values were chosen through trial and error.

A prediction window length of 1 second was used with 16 coincidence points.  A 1 second window was deemed sufficient for the purposes of angular rate following as it is assumed that the angular rate demands are constant across the window length.  This is a reasonable assumption as the angular rates do not change significantly after 1 second.

\subsection{Fault Simulation}
The concept behind the fault tolerant controller design for the 6DoF model is based on monitoring the control derivatives.  The non-dimensional aerodynamic coefficients for the forces and moments given in appendix I are made up of a series of aerodynamic and control derivatives.  For example the term $- 6.54\times 10^{-3}\delta_e$ in the pitching moment coefficient represents the pitch control derivative, $C_{m_{\delta_e}}$, the contribution of the elevator control input on the pitching moment coefficient.  In the example given above the value $- 6.54\times 10^{-3}$ is specific to the given aircraft as are all the derivative values given in appendix I.  For any aircraft these values are obtained via experimental testing or computational fluid dynamic techniques, and as the derivatives are affected by any physical change in the control surface any change in a control derivative would indicate a fault. 

For the simulation results given in subsection \ref{subsection:chap5_6DoFNR} the faults are simulated by reducing the efficiency of the control surface.  The primary role of an elevator is to provide pitch control, so its largest contribution is on the pitching moment, and therefore a change in the $C_{m_{\delta_e}}$ derivative would indicate an elevator fault.  The aileron contributes primarily to the rolling moment $C_l$ and the control derivative associated with the aileron from equation \eqref{eqn:chap5_Cl} is $C_{l_{\delta_a}} = 6.1\times 10^{-4}$.  Finally the rudder has the biggest impact on the yawing moment, $C_n$, and the associated control derivative from equation \eqref{eqn:chap5_Cn} is $C_{n_{\delta_r}} = - 9.0\times 10^{-4}$.  To simulate a fault in a control surface the respective control derivative is reduced. 

\subsection{Numerical Results}\label{subsection:chap5_6DoFNR}
\begin{figure}[H]
\hspace{-0.3in}
\includegraphics[scale=0.3]{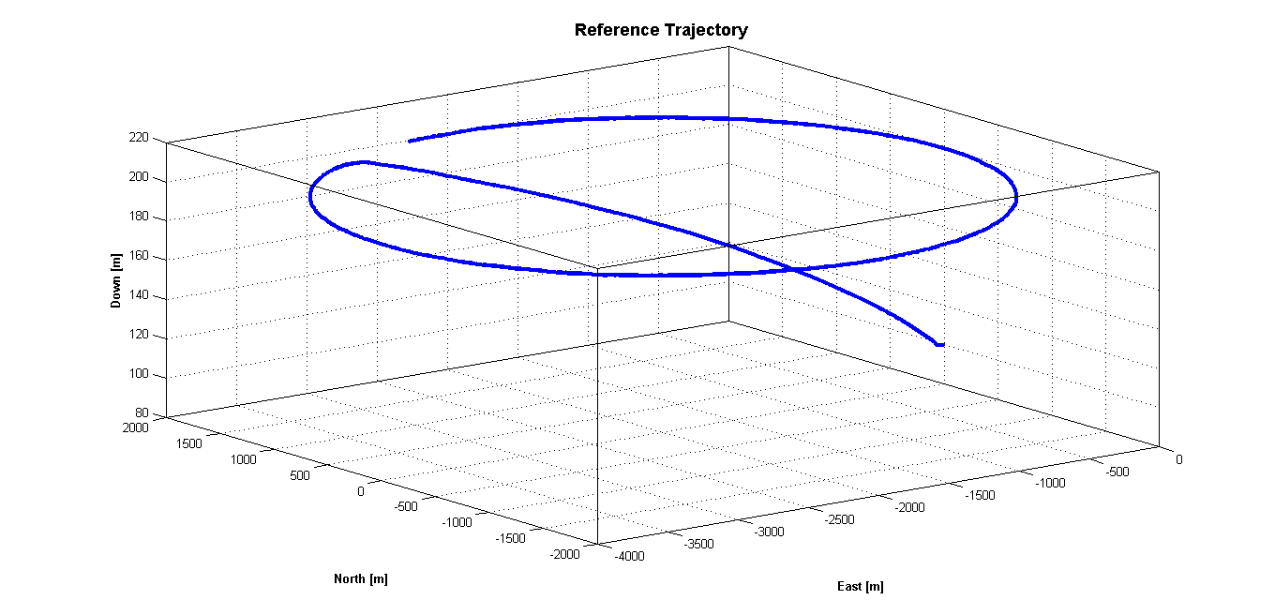} 
\caption{6DoF Reference Trajectory}
\label{fig:chap5_refTraj_6DoF}
\end{figure}

To investigate the effectiveness of the NMPC controller design as a fault tolerant controller the aircraft was required to fly the trajectory given in figure \ref{fig:chap5_refTraj_6DoF}. 

Three different scenarios were set up:

\begin{enumerate}[label=\bfseries Scenario \arabic*:, leftmargin = 100pt]
\item faulty elevator: $70\%$ reduction in efficiency 20 seconds into flight,
\vspace{-10pt}
\item faulty aileron: $80\%$ reduction in efficiency 20 seconds into flight,
\vspace{-10pt}
\item faulty rudder: $60\%$ reduction in efficiency 20 seconds into flight,
\end{enumerate}

Each scenario is run with and without FDI information.  The FDI information when used is assumed and has not been modelled.  It is assumed that the FDI subsystem is capable of providing the time of fault and the efficiency of the control surface.

\subsection{Analysis}
\begin{figure}[H]
\hspace{-0.4in}
%\hspace{-0.7in}
\includegraphics[scale=0.4]{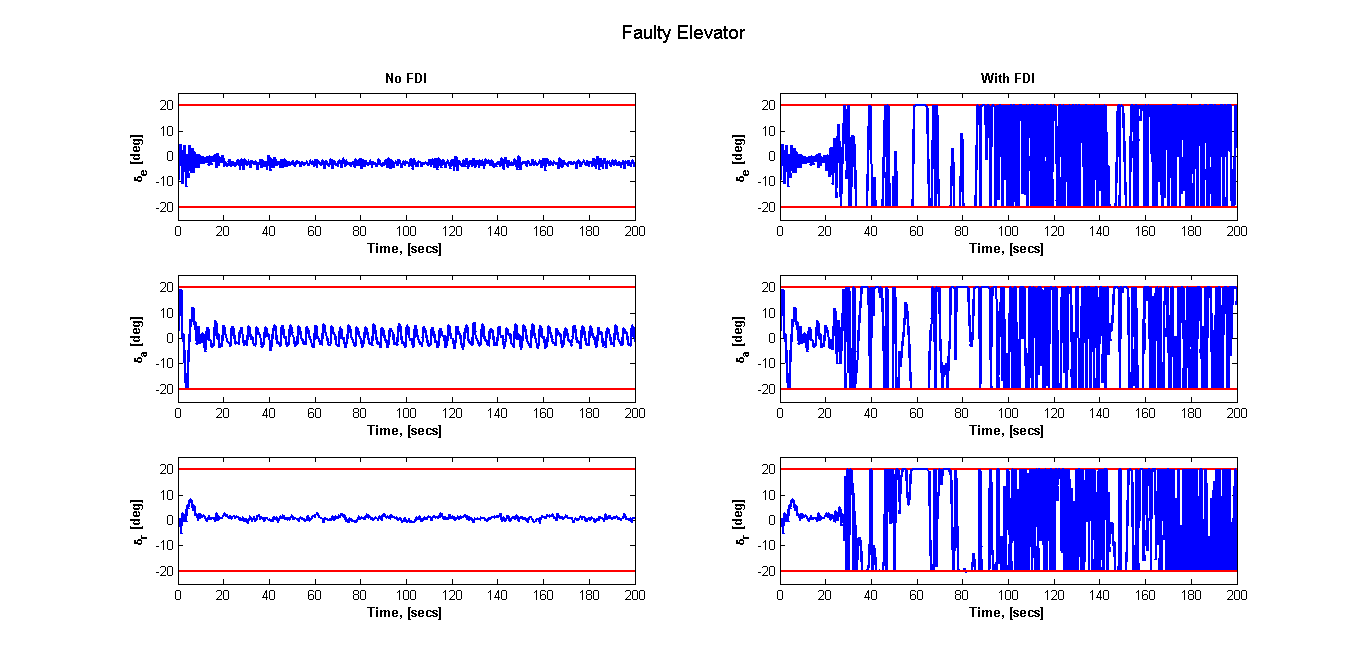} 
%\vspace{-0.5in}
\caption{Faulty Elevator: Control Surface Activity, constraints (red), control surface activity (blue).  Left column: no FDI information, Right column: with FDI information}
\label{fig:chap5_FE}
\end{figure}

\begin{figure}[H]
\hspace{-0.4in}
%\hspace{-0.7in}
\includegraphics[scale=0.4]{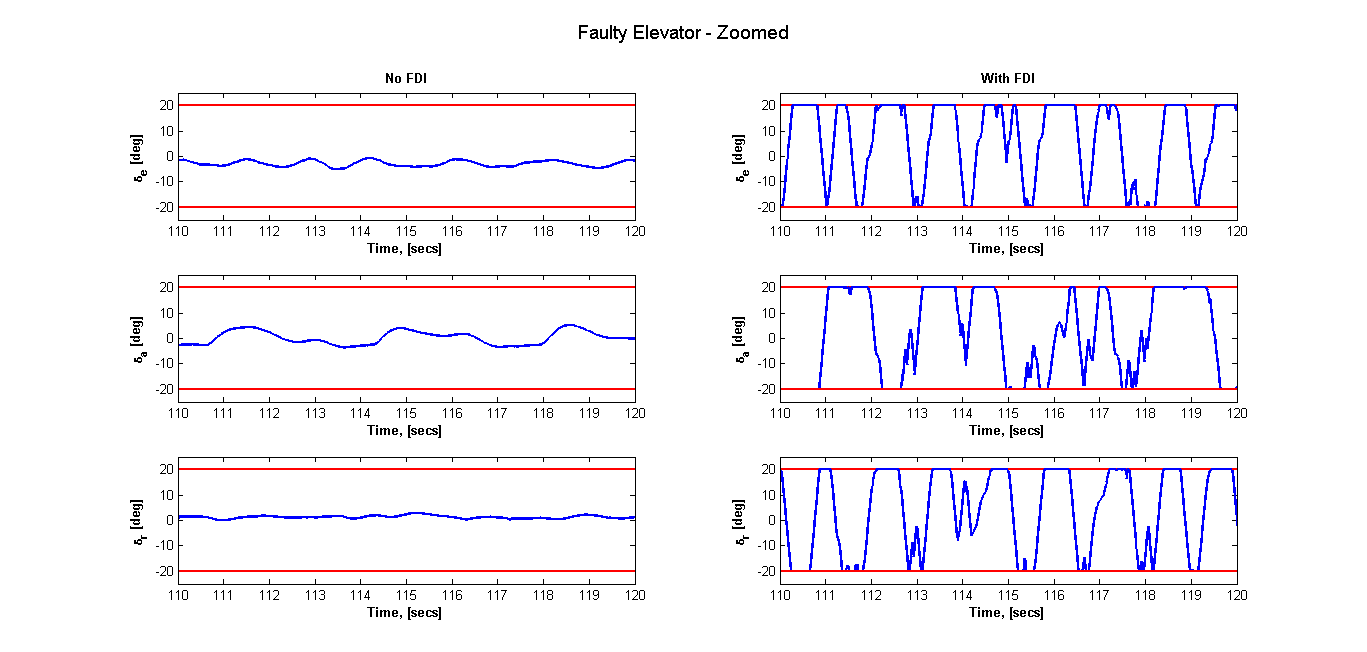} 
%\vspace{-0.5in}
\caption{Faulty Elevator: Control Surface Activity (zoomed), constraints (red), control surface activity (blue).  Left column: no FDI information, Right column: with FDI information}
\label{fig:chap5_FE_zoomed}
\end{figure}

Figure \ref{fig:chap5_FE} presents the plots for the control surface activity given an elevator with $70\%$ reduction in efficiency.  Figure \ref{fig:chap5_FE_zoomed} presents a 10 second plot to show a close up of the elevator activity.  The plots show that without any knowledge of the fault the activity in the elevator decreases after 20 seconds and there is very little change in the aileron and rudder activity once the fault occurs.  When FDI information is provided however the knowledge of the fault prompts the control surfaces to work harder to compensate for the fault.  This is seen in all three control surfaces which at various times during the flight are all operating at the constraints. 

\begin{figure}[H]
\hspace{-0.4in}
%\hspace{-0.7in}
\includegraphics[scale=0.4]{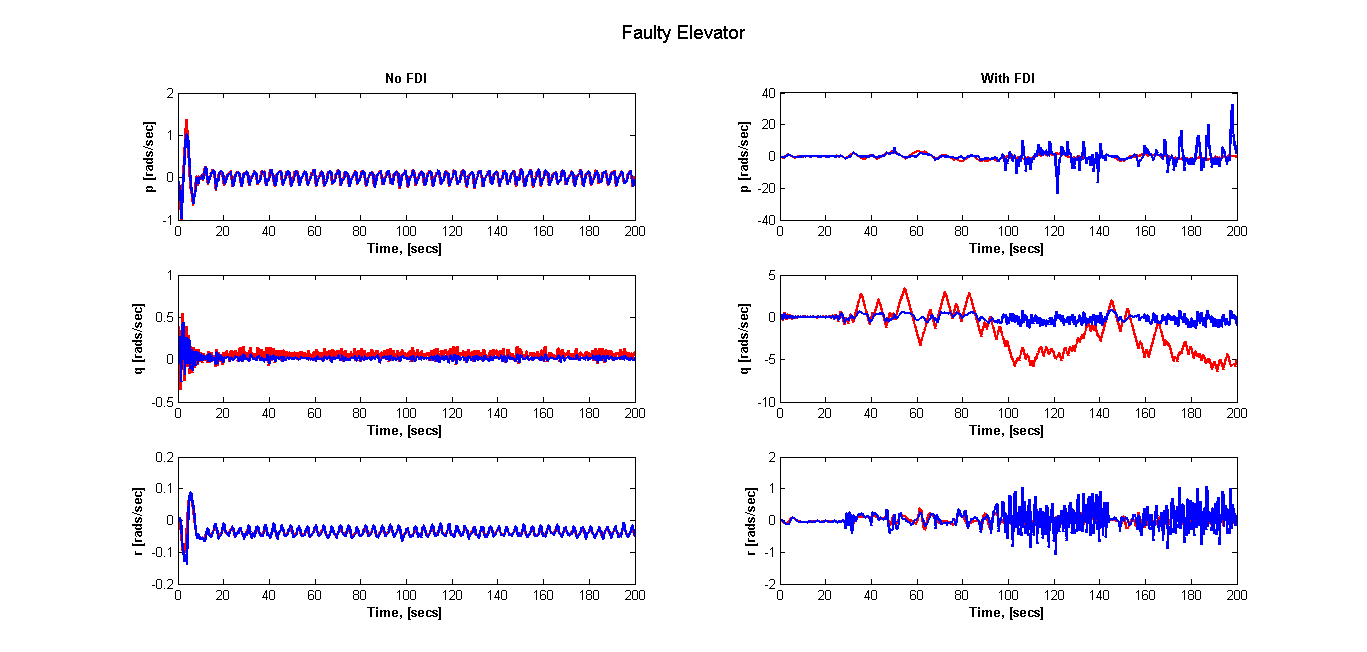}
%\vspace{-0.5in}
\caption{Faulty Elevator: Angular Rates, demanded (red), actual (blue).  Left column: no FDI information, Right column: with FDI information}
\label{fig:chap5_AR_FE}
\end{figure}

The angular rate plots are shown in figure \ref{fig:chap5_AR_FE}.  A fault in the elevator directly affects the pitch rate $q$, and without any FDI information the controller is unable to meet the pitch rate demands, however the roll rate and yaw rate demands are followed very closely.  With knowledge of the faults there is an increase in the demanded angular rates and the controller shows a significant improvement in performance in being able to follow the demanded rates.

\begin{figure}[H]
\hspace{-0.3in}
%\hspace{-1in}
\includegraphics[scale=0.4]{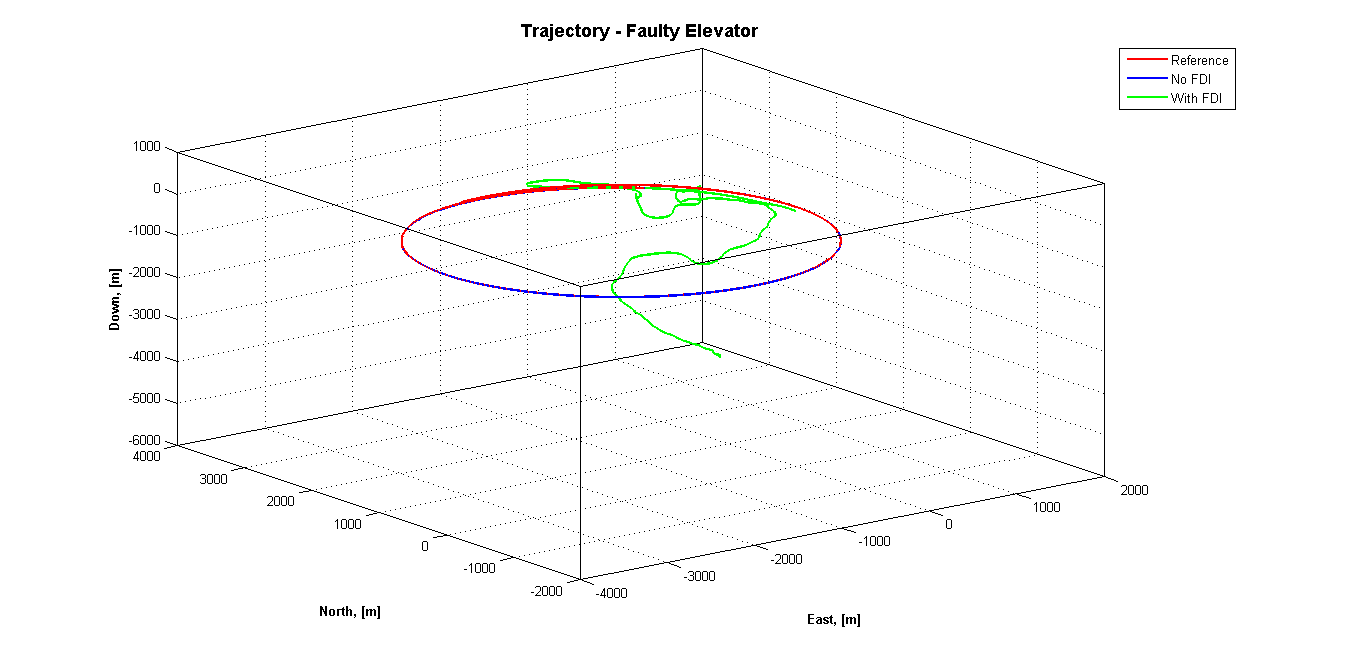} 
%\vspace{-0.5in}
\caption{Faulty Elevator: 6DoF Trajectory}
\label{fig:chap5_traj_FE}
\end{figure}

The trajectories flown by the aircraft with a faulty elevator, with and without FDI information, are provided in figure \ref{fig:chap5_traj_FE}.  The results show that in the absence of FDI information the aircraft successfully flies the trajectory, however providing FDI information caused the simulated trajectory to diverge from the reference trajectory.  This result shows that the controller behaves exactly as expected.  The controller has been designed to maintain the angular rate demands not the reference trajectory.  The angular rate plots show that with the FDI information there is an increase in performance of the controller in terms of tracking the angular rate demands.  The trajectory plots show that the simulated trajectory produced with FDI information causes the aircraft to drop below ground level which is physically impossible.  This is a result of not applying constraints on the aircraft position vector.  Hence unless a parameter is explicitly penalised in the cost function and/or constraints placed upon the parameters the controller will use everything available to it to achieve what is being demanded of it.

The same plots were produced for scenario 2 but have been omitted due to space constraints.  The aileron and rudder are primarily used to control the lateral motion of the aircraft while the elevator controls the longitudinal motion.  Thus the results showed very little change in the behaviour of the elevator when FDI information is provided compared to no FDI information.  The rudder and aileron on the other hand when provided with FDI information increase their activity after the occurrence of the fault to compensate for the loss in efficiency and operate closer to the constraints.    

Similar behaviour was present in the angular rate plots.  There was little or no change in the pitch response of the aircraft once FDI information is provided compared when FDI is absent.  In the case of no FDI the actual roll rate is lower than the demand, however once information on the fault is provided tracking performance increases.  This is also true for the yaw rate response.  In the presence of an aileron fault roll and yaw rate demands increase to sustain lateral motion. 

The trajectory plots showed that in the case where FDI information is provided the aircraft was seen to deviate slightly off the path.  The deviation is not as significant in the event of an aileron fault as the rudder also helps to control the lateral motion of the aircraft hence providing an extra degree of redundancy.  

The same plots were produced for scenario 3 with a $60\%$ reduction in efficiency in the rudder and it was again observed that as the elevator has very little influence on lateral motion, there was little change in elevator activity with no difference between the no FDI and with FDI cases.  The rudder was seen to be pushed to its lower limit and the aileron increased in the negative direction causing the aircraft to bank more to the left.  

A faulty rudder was seen to have no effect on the angular rate demands.  Tracking performance was the same both with and without FDI information.  This was translated in the trajectory plots which illustrated that the aircraft closely followed the flight path with and without FDI information.  

Another point to note is that the non-dimensional force and moment coefficients given in appendix I are only valid for angle of attack $\alpha \leq 15^0$.  Plots of $\alpha$ were produced for all scenarios and the results showed that $\alpha$ never went above $15^0$.  Models of $\alpha > 15^0$ are also given in \cite{garza2003collection}, so to cover all flight envelopes multiple $\alpha$ models could have been incorporated or constraints could have been placed on $\alpha$ to ensure $15^0$ is never exceeded.  In our work, we added stall speed constraints to the NMPC controller and the results show that stall speed was never encountered.

Overall the results of our 6DoF analysis show that NMPC design as a fault tolerant controller is viable, showing that in the absence of FDI information the controller is capable of allocating control authority to the appropriate actuators to fly the aircraft on the given flight path.  This is a display of the inherent fault tolerant capabilities of NMPC.  Turning on FDI updates improved the tracking performance of the controller.  The results did however show that unless a quantity is penalised in the cost function, and/or constraints are applied, the controller will push the limits to achieve the desired outcome.  In this case the controller was specifically designed to track angular rate demands, hence providing FDI information resulted in an increase in tracking performance in the event of a control surface fault.

The next section will look at the design of an FDI filter to be incorporated into the FTC developed in this section. 

%******************End of old chapter 10*****************************

%*****************Begining of old chapter 11************************
\section{Fault Detection and Identification}\label{section:FDI}
The fault detection concepts from \cite{RKPaper2} are implemented here for the full 6DoF aircraft model and designs for the UKF are presented.  To design the filter a PID controller for the aircraft was designed and implemented before integrating the filter with the NMPC controller.

A traditional PID controller was used to control the aircraft through the range of manoeuvres required to test and tune the filter.  The PID control method, although not optimal in terms of performance, was quick to implement and tune to the level required.

The proposed fault detection scheme is based on the principle that a failure in any one of the control surfaces would directly affect the corresponding control derivative.  Hence changes in the control derivatives would indicate a fault has occurred, while at the same time the filter would provide the controller with estimates of the derivatives.  Furthermore, up to date estimates of the derivatives will allow the NMPC controller to perform at its optimum.

The force and moment equations given in appendix I show that there are a total of 24 control derivative.  These are listed in Table \ref{table:chap5_CDs}. 

\begin{table}[H]
\caption{Control Derivatives}
\label{table:chap5_CDs}
\begin{center}
\begin{tabular}{|c|c||c|c|}
\hline 
\textbf{Derivative} & \textbf{Value} & \textbf{Derivative} & \textbf{Value} \\ 
\hline 
$CX_{dE1}$ & $9.5\times 10^{-4}$ & $Cl_{dR2}$ & $4.5\times 10^{-6}$  \\ 
\hline 
$CX_{dE2}$ & $8.5\times 10^{-7}$ & $Cl_{dE1}$ & $5.24\times 10^{-5}$ \\ 
\hline 
$CY_{dE1}$ & $1.75\times 10^{-4}$ & $Cm_{dE1}$ & $6.54\times 10^{-3}$ \\ 
\hline 
$CY_{dR1}$ & $1.55\times 10^{-3}$ & $Cm_{dE2}$  & $8.49\times 10^{-5}$ \\ 
\hline 
$CY_{dR2}$ & $8\times 10^{-6}$ & $Cm_{dE3}$ & $3.74\times 10^{-6}$ \\ 
\hline 
$CZ_{dE1}$ & $4.76\times 10^{-3}$ & $Cm_{dA1}$  & $3.5\times 10^{-5}$  \\ 
\hline 
$CZ_{dE2}$ & $3.3\times 10^{-5}$ & $Cn_{dA1}$  & $1.4\times 10^{-5}$ \\ 
\hline 
$CZ_{dA1}$ & $7.5\times 10^{-5}$ & $Cn_{dA2}$  & $7.0\times 10^{-6}$ \\ 
\hline 
$Cl_{dA1}$ & $6.1\times 10^{-4}$ & $Cn_{dE1}$  & $8.73\times 10^{-5}$  \\ 
\hline 
$Cl_{dA2}$ & $2.5\times 10^{-5}$ & $Cn_{dE2}$  & $8.7\times 10^{-6}$  \\ 
\hline 
$Cl_{dA3}$ & $2.6\times 10^{-6}$ & $Cn_{dR1}$  & $9.0\times 10^{-4}$ \\ 
\hline 
$Cl_{dR1}$ & $-2.3\times 10^{-4}$ & $Cn_{dR2}$ & $4.0\times 10^{-6}$\\
\hline
\end{tabular} 
\end{center}
\end{table}

\begin{figure}[H]
\hspace{-0.4in}
\includegraphics[scale=0.3]{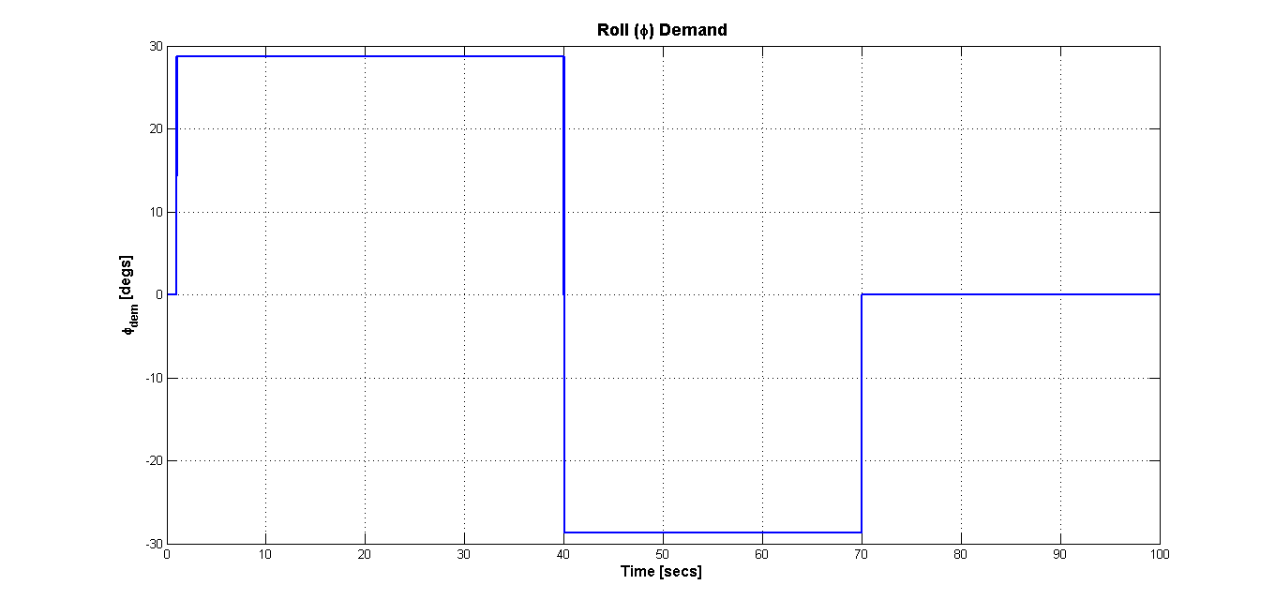} 
%\vspace{-0.5in}
\caption{6DoF Motion Filter Tests - Roll Angle Demands}
\label{fig:chap5_phiDem}
\end{figure}

To test the filters the aircraft was required to achieve the roll angle demands given in figure \ref{fig:chap5_phiDem}. 

\begin{figure}[H]
\hspace{-1in}
\includegraphics[scale=0.4]{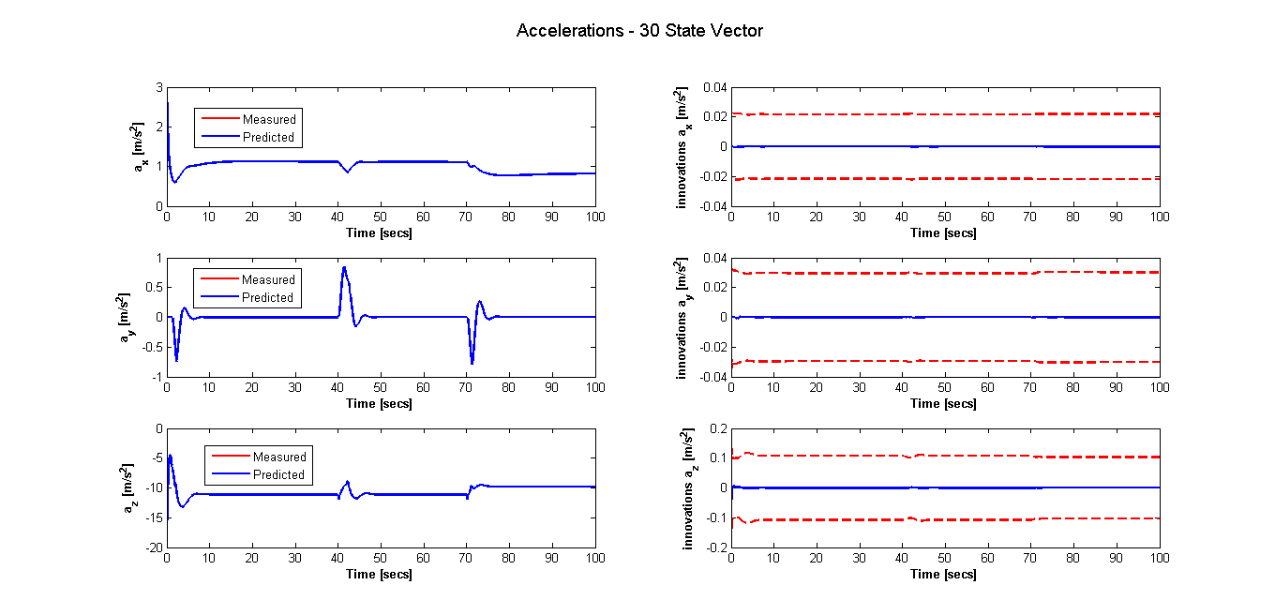}
%\vspace{-0.6in}
\caption{Accelerations - 30 State Vector, left column estimates (red measured, blue predicted), right column innovations (2$\sigma$ uncertainty bounds (red dashed) and innovations (blue).}
\label{fig:chap5_Accel_fullFilter}
\end{figure}

\begin{figure}[H]
\hspace{-1in}
\includegraphics[scale=0.4]{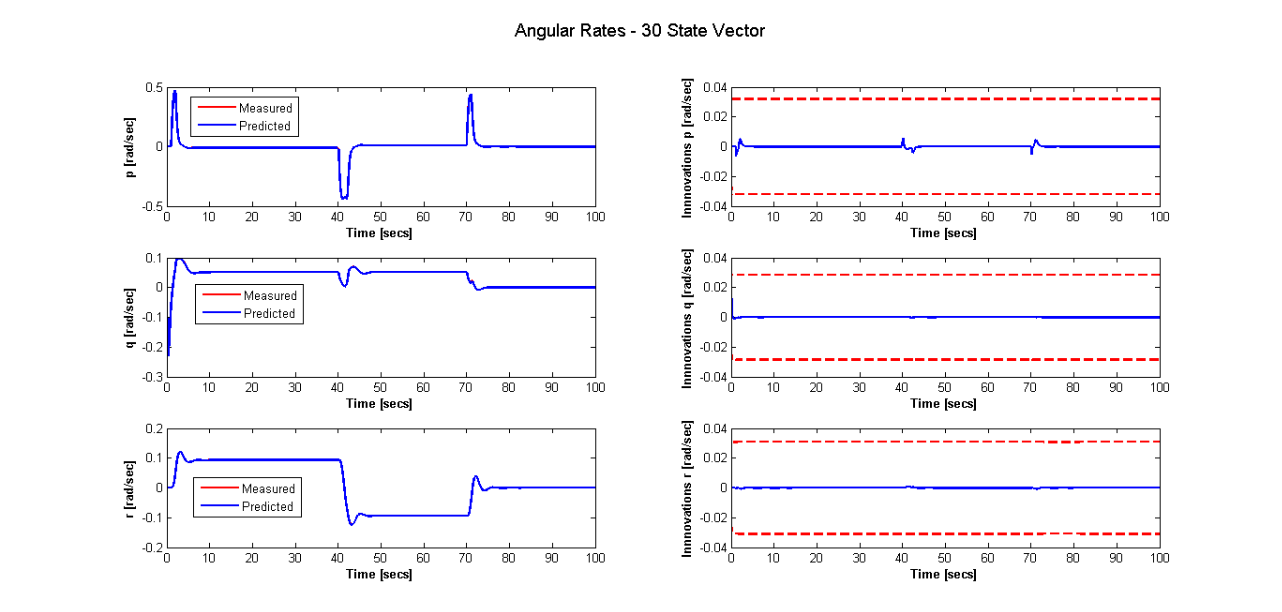}
%\vspace{-0.6in}
\caption{Angular Rates - 30 State Vector, left column estimates (red measured, blue predicted), right column innovations (2$\sigma$ uncertainty bounds (red dashed) and innovations (blue).}
\label{fig:chap5_AR_fullFilter}
\end{figure}

\begin{figure}[H]
\hspace{-1in}
\includegraphics[scale=0.4]{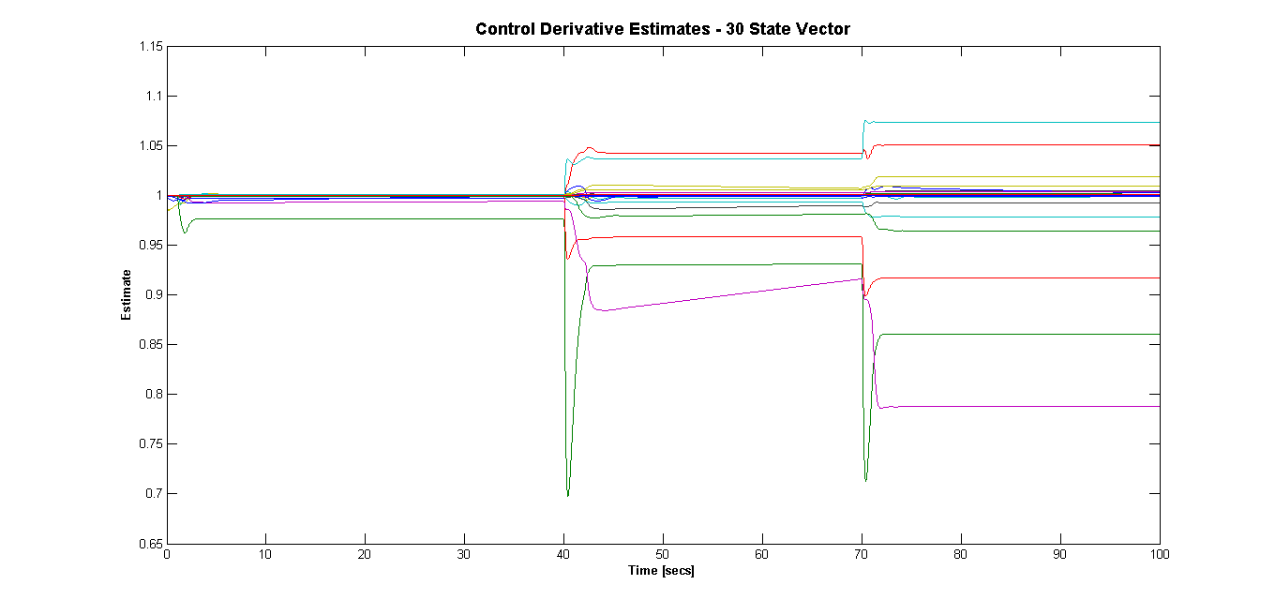}
%\vspace{-0.6in}
\caption{Control Derivative Estimates - 30 State Vector.  Each line corresponds to a normalised value of a control derivative estimate and should have a value of 1.}
\label{fig:chap5_CDests_fullFilter}
\end{figure}

Initially a 30 state UKF filter was designed where the states comprised of three accelerations $\left(a_x,\,a_y,\,a_z\right)$, three angular rates $\left(p,\,q,\,r\right)$ and the 24 control derivatives given above.  The measurements were of the body acceleration and angular rates (as would be provided by an IMU sensor).  All the derivatives were normalised to 1 hence the states of the control derivative were set to 1.  The results of the acceleration and angular rate innovations are given in figures \ref{fig:chap5_Accel_fullFilter} and \ref{fig:chap5_AR_fullFilter} respectively.  The results show that the filter does an excellent job of predicting the accelerations and angular rates as the filter predictions align perfectly with the measurements of acceleration and angular rates.  The estimates of the control derivatives are shown in figure \ref{fig:chap5_CDests_fullFilter}.  Since all derivatives were normalised the estimates should all have a value of 1.  However, as the plot shows, the filter was unable to correctly estimate the value of the derivatives, as many of the states in the filter are unobservable.


To address the issue of observability the number of states was reduced to 19; 3 accelerations, 3 angular rates and 13 control derivatives.  This was achieved by having one control derivative estimate per force/moment for a particular control surface.  For example, the control derivatives given in table \ref{table:chap5_CDs} show that there are 2 $CX$ control derivatives which are due to the elevator $CX_{dE1}$ and $CX_{dE2}$; the 19 state filter has only one derivative $CX_{dE}$ used to represent both of the $CX$ derivatives due to the elevator.  Thus the contributions of each control surface in the force and moment equations are grouped together in this manner reducing the number of control derivative states from 24 to 13.  The plots of angular rates and accelerations have been omitted due to space constraints, with the results again showing close to perfect compliance between prediction and measurement.  However while plots of the control derivative estimates did show an improvement in estimates the issue of unobservable states was still evident.    

To further address this issue the number of filter states was again reduced, from 19 to 12.  For control surface failure the most important derivatives are deemed to be $Cl_{dA1}$ for aileron, $Cm_{dE1}$ for elevator and $Cn_{dR1}$ for the rudder, and the trim values $Cl_0$, $Cm_0$ and $Cn_0$.  From the equations given in Appendix I the trim values corresponding to aileron, elevator and rudder are $Cl_0 = 0$, $Cm_0 = -6.61\times 10^{-3}$ and $Cn_0 = 0$ respectively.

The plots obtained of acceleration and angular rate estimates for the 12 state filter showed that the angular rate estimations were excellent however, the filter was unable to make correct estimates of acceleration.  This was to be expected because the filter states are unrelated to force, and are all moment related terms.  The control derivatives $Cl_{dA1}$, $Cm_{dE1}$ and $Cn_{dR1}$ are normalised to 1 as is the trim value for $Cm_0$ and the estimates of these were plotted.  Results showed big discrepancies between the actual and estimated values for the elevator terms.          

Due to the continuing presence of unobservable states the state vector was once more reduced by removing the acceleration terms resulting in a 9 state filter.  The measurements supplied to the filter were only of the angular rates.  The acceleration terms were removed due to the errors present in the estimates.  The angular rate estimates and innovations were plotted and as expected, show that the filter predictions closely match the measurements.  It was evident however from the control derivative and trim estimate plot that the observability issue was still present.     

In a final attempt to solve the observability issue three separate filters were developed, one each for roll, pitch and yaw with each filter, a 3 state filter.  The states for the roll only filter are $p,\,Cl_{dA1}$ and $Cn_0$, the pitch only filter states are $q,\,Cm_{dE1}$ and $Cm_0$ and the yaw only filter has $r,\,Cn_{dR1}$ and $Cn_0$ as states.  The angular rate estimates were obtained and again showed compliance with the measurements.  Separating the filters caused slight improvement in the control derivative and trim estimates however the observability problem was still present particularly for the $Cm_0$ term.  This was to be expected as the aircraft lateral dynamics have been excited by the demanded roll inputs given in figure \ref{fig:chap5_phiDem} hence the estimates of the derivatives related to the later dynamics are more accurate than the longitudinal motion derivatives.  For good estimates it is necessary to excite the aircraft dynamics. 

The results obtained in the previous section illustrate that if the NMPC controller could be provided with estimates of the control derivatives it would assist the controller in allocating control authority appropriately.  For this reason a UKF filter was designed in this section to provide real time estimates.  However results showed that many of the control derivatives are unobservable.  Many attempts were made to tackle this issue, however all proved to be unsuccessful.  The 3 filter solution was integrated with the NMPC controller to test the full active fault tolerant control system.  Results for these tests have not been supplied as the incorrect estimates of the filter caused the solution from the controller to diverge.  Further investigations into the full 6DoF fault tolerant controller are required.

As the main objective of this work is to demonstrate fault tolerant flight control, the next section looks at the longitudinal motion of the aircraft with integrated FDI to form a full active fault tolerant flight control system.    

%*********************End of old chapter 11***********************
\section{Engine Failure - Loss of Power} \label{section:engineFailure}
The loss of power on an aircraft due to engine failure can result in a catastrophic breakdown of the system if left unattended.  This section demonstrates the use of the active FTC system design as a fault tolerant flight controller in the event of an engine failure.

The FTC system is used to control the longitudinal motion of the aircraft.  The design comprises of a UKF filter to monitor the thrust level of the air vehicle.  Fault detection logic is built into the filter and once the decision is made that there is a loss of power the filter estimates are fed to the NMPC controller for reconfiguration.

The filter design is detailed in the next subsection followed by the details of the controller design in the subsection after.  Finally numerical results are presented.

\subsection{Filter Design}
The filter design process consists of the development of a simple PID thrust controller.  The NMPC and filter designs were independently constructed and tested then integrated into the final design.

\subsubsection{UKF FDI Filter}
The UKF filter is designed to estimate the amount of thrust used by the aircraft.  The filter states are:

\begin{equation}
\mathbf{x_{\text{ukf}}} = [V_N,\,\,V_D,\,\,\theta,\,\,T],
\end{equation}    

where $T$ is thrust (see appendix II for the thrust model).  The measurement vector is:

\begin{equation}
\mathbf{x_{\text{ukf}}} = [V_\text{EAS},\,\,v_D,\,\,\theta],
\end{equation}   

where $V_\text{EAS}$ is equivalent airspeed of the aircraft at sea level whereas $V_T$ is the true airspeed at altitude.  For this work the aircraft is assumed to be flying low enough for $V_\text{EAS} = V_T$.

The weighting matrices $Q$ and $R$ were set to:

\begin{equation}
\mathbf{Q} = \begin{bmatrix}
(5\,\Delta t\,0.05)^2 & 0 & 0 & 0\\
0 & (5\,\Delta t\,0.05)^2 & 0 & 0\\
0 & 0 & (0.1\,\Delta \, t)^2 & 0\\
0 & 0 & 0 & (6500\,\Delta t \,0.3)^2
\end{bmatrix},
\end{equation}
\newline
\newline
\begin{equation}
\mathbf{R} = \begin{bmatrix}
(0.05)^2 & 0 & 0\\
0 & (0.05)^2 & 0\\
0 & 0 & (0.017)^2
\end{bmatrix},
\end{equation}
\newline
where $\Delta t$ is the filter update rate 0.01 secs.  The initial state vector and covariance matrix are:

\begin{equation}
\mathbf{x}(0) = \left[50,\,\,0\,\,0.04247,\,\,1507.7526\right]^\intercal, \quad
\mathbf{P}(0) = \begin{bmatrix}
(0.5)^2 & 0 & 0 & 0\\
0 & (0.5)^2 & 0 & 0\\
0 & 0 & (0.017)^2 & 0\\
0 & 0 & 0 & (315)^2
\end{bmatrix}
\end{equation}

\subsubsection{Numerical Results}  
The following test cases were carried out to examine the filter performance:

\begin{enumerate}[label=\bfseries Test \arabic*:, leftmargin = 100pt]
\item no fault,
\vspace{-10pt}
\item $70\%$ loss of power 70 secs into flight,
\vspace{-10pt}
\item $90\%$ loss of power 35 secs into flight,
\vspace{-10pt}
\item $50\%$ loss of power 20 secs into flight.
\end{enumerate}

\begin{figure}[H]
\hspace{-0.4in}
\includegraphics[scale=0.3]{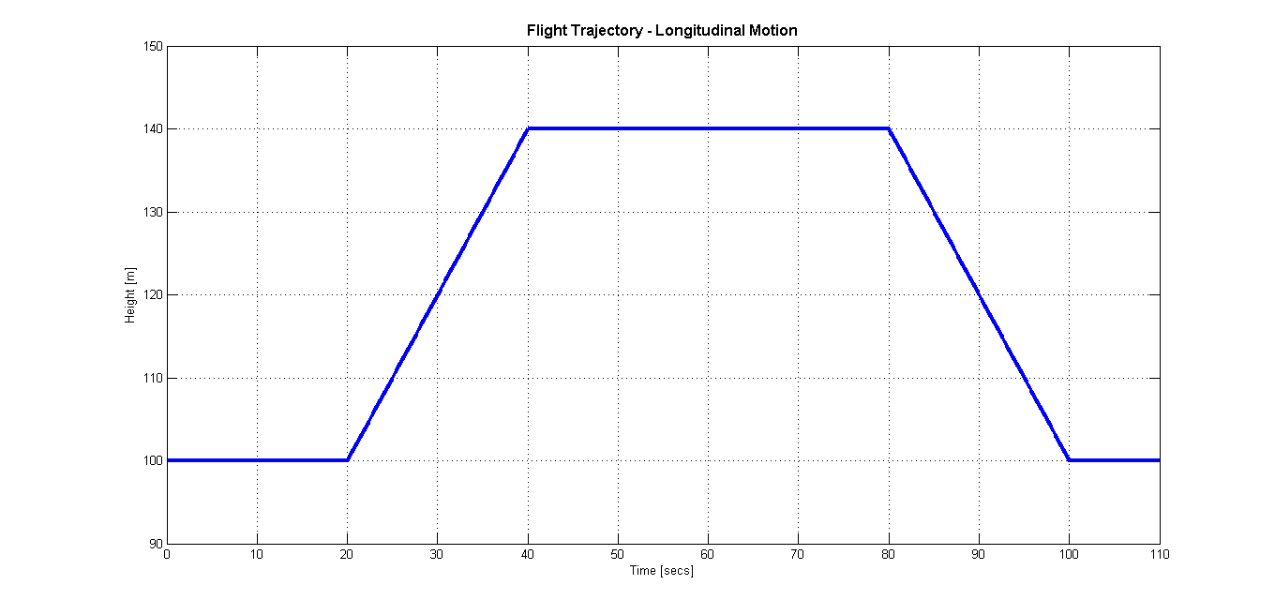}
\caption{Reference trajectory for longitudinal motion}
\label{fig:long_traj} 
\end{figure}

The aircraft was required to fly the trajectory given in figure \ref{fig:long_traj}.  The effects of wind and turbulence have been taken into account as well as the effect of noise on the measurements of $V_{\text{EAS}}$, $v_d$ and $\theta$, modelled as a normally distributed random white noise.

\begin{figure}[H]
\hspace{-1in}
\includegraphics[scale=0.4]{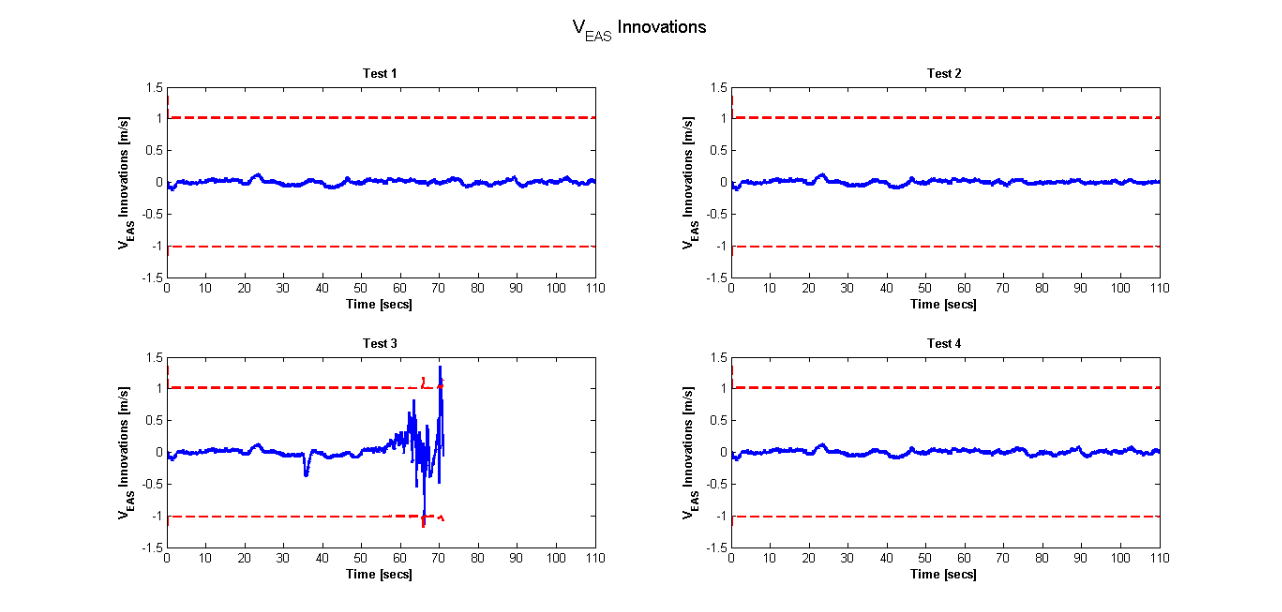} 
%\vspace{-0.55in}
\caption{UKF $V_{\text{EAS}}$ Innovations - Longitudinal Model, $\pm 2\sigma$ innovation covariance bounds (red dashed lines), $V_{\text{EAS}}$ innovations (solid blue line)}
\label{fig:chap5_VEASInnovations_longThrust}
\end{figure}

\begin{figure}[H]
\hspace{-1in}
\includegraphics[scale=0.4]{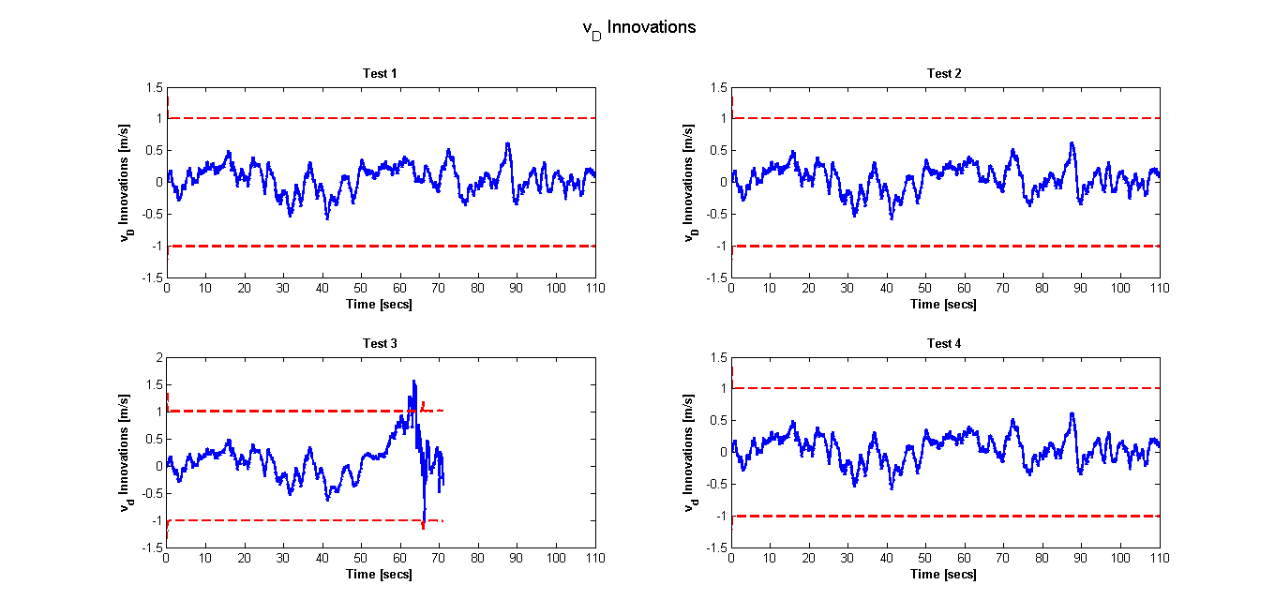} 
%\vspace{-0.55in}
\caption{UKF $V_D$ Innovations - Longitudinal Model, $\pm 2\sigma$ innovation covariance bounds (red dashed lines), $V_D$ innovations (solid blue line)}
\label{fig:chap5_vDInnovations_longThrust}
\end{figure}

\begin{figure}[H]
\hspace{-1in}
\includegraphics[scale=0.4]{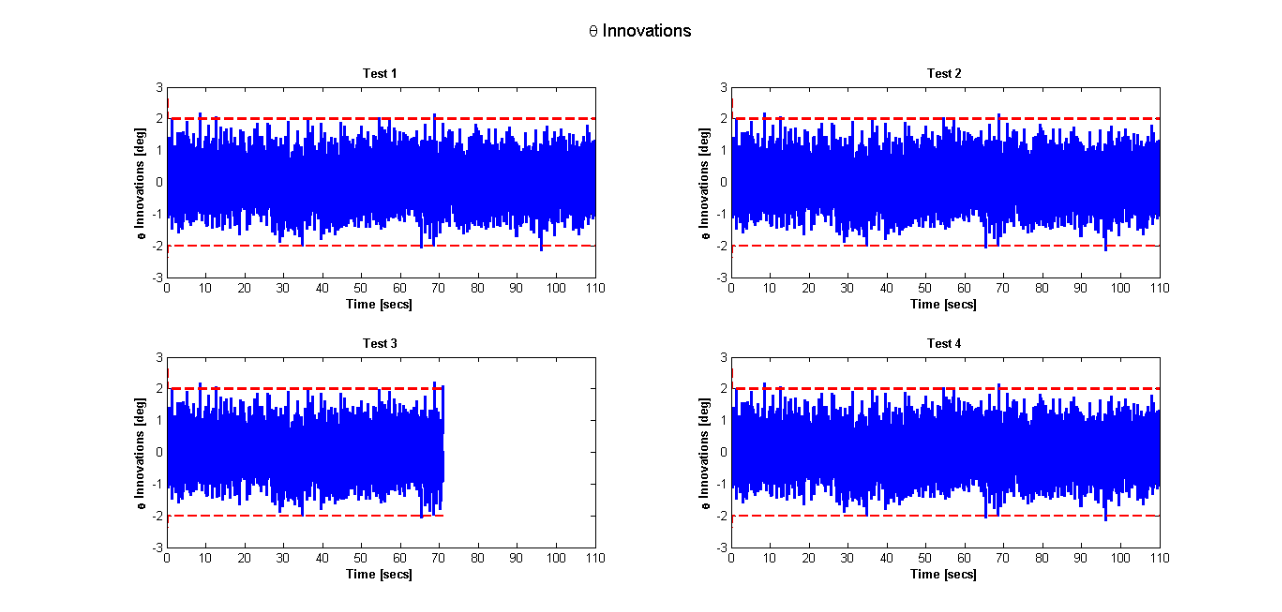} 
%\vspace{-0.55in}
\caption{UKF $\theta$ Innovations - Longitudinal Model, $\pm 2\sigma$ innovation covariance bounds (red dashed lines), $\theta$ innovations (solid blue line)}
\label{fig:chap5_thetaInnovations_longThrust}
\end{figure}

\begin{figure}[H]
\hspace{-1in}
\includegraphics[scale=0.4]{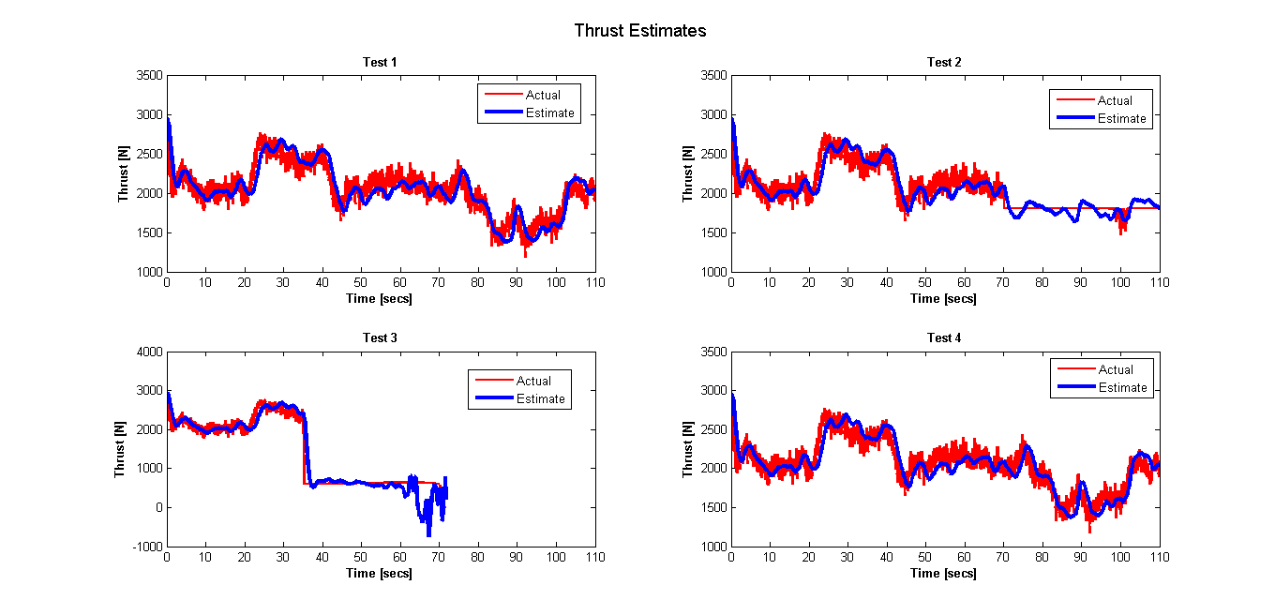}
%\vspace{-0.55in} 
\caption{UKF Thrust Estimates - Longitudinal Model}
\label{fig:chap5_thrustEstimates}
\end{figure}

To analyse the performance of the filter the innovation covariance plots were examined and are given in figures \ref{fig:chap5_VEASInnovations_longThrust}, \ref{fig:chap5_vDInnovations_longThrust} and \ref{fig:chap5_thetaInnovations_longThrust} for  $V_{\text{EAS}}$, $v_d$ and $\theta$ respectively.  The results show that for all test cases the innovations are well within the $2\sigma$ covariance bounds.  The test case 3 where the thrust level drops to $10\%$ shows that after approximately 70 seconds the aircraft is unable to maintain flight as there is not enough power hence the filter diverges.  The thrust estimates are given in figure \ref{fig:chap5_thrustEstimates} along with the actual thrust applied to the aircraft.  In each test case the filter does an excellent job of estimating the thrust levels.  


\subsubsection{Fault Detection Logic}\label{subsubsec:chap5_FDL}
The premise behind our FTC design is the provision of the updates of the power status to the NMPC controller to enable controller reconfiguration.  When an engine fails the amount of thrust available decreases.  If this level of thrust could be estimated and provided to the NMPC controller the maximum constraint on thrust can be updated and the controller can then allocate control authority to the control inputs accordingly.  For this reason it is important to detect the fault and to know when to begin feeding the controller with filter estimates of thrust, hence the need for fault detection logic.

The controller is designed (section \ref{subsec:chap5_thrust_NMPC})  to calculate the optimal amount of thrust to maintain a height demand and true airspeed.  The filter on the other hand estimates the thrust level currently used by the aircraft, hence if the demand is greater than the estimate this would indicate an engine failure.

\begin{figure}[H]
\hspace{-1in}
\includegraphics[scale=0.4]{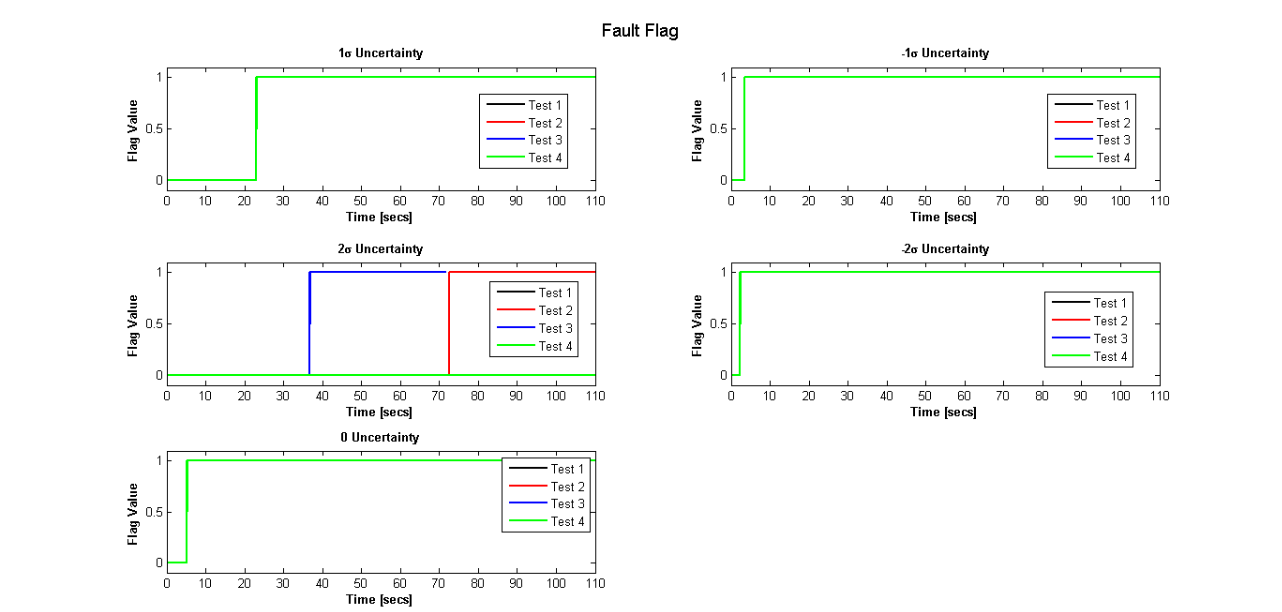}
%\vspace{-0.4in} 
\caption{UKF Fault Flag - Longitudinal Model}
\label{fig:chap5_faultFlag}
\end{figure}

The fault detection logic therefore comprises of checking whether the thrust demand is higher than the thrust estimate.  If this is true for a set period of time then a fault has occurred and a flag is turned on indicating that a fault has occurred and consequently the constraints must be updated via the filter estimates.  The filter outputs an estimate of the state as well as the uncertainty on the estimate, so that the actual value of the state as predicted by the filter is within plus or minus the uncertainty.  For this reason a number of tests were performed to see whether the check should include zero level of uncertainty, $\pm 1\sigma$ uncertainty or $\pm 2\sigma$, with the results given in figure \ref{fig:chap5_faultFlag}.  The results are based on a fault count of 200, i.e. when the demand is greater than the thrust estimate the fault counter is incremented by one, and when this counter exceeds 200 the fault flag switches from 0 to 1 indicating to the controller that the maximum constraint on thrust must be updated with the filter estimate.  The number of counts being set to 200 is based purely on trial and error.  The results show that the filter estimate plus $2\sigma$ uncertainty was able to correctly identify the fault within approximately a couple of seconds of the fault occurring.  The other uncertainty bounds as well as the zero uncertainty case all indicated false detection of the fault, that is the fault flag is set to true at the incorrect times.  Note that a fault was not detected for test case 4 even with a $2\sigma$ uncertainty bound.  This is because the thrust estimate plots (as shown in figure \ref{fig:chap5_thrustEstimates}) indicate that in a no fault case the aircraft requires no more than $50\%$ of the maximum thrust to maintain the given trajectory, hence the demand is at all times less than the estimate.

In the next subsection the complete active FTC system design for thrust control is detailed.

\subsection{Controller Design}\label{subsec:chap5_thrust_NMPC}
This section steps through the controller design for the active FTC system for the longitudinal motion of the aircraft.

Pseudospectral discretisation \cite{RKPaper1} is applied to the controller design and the NMPC state vector is:

\begin{equation}\label{eqn:xnmpc_thrust}
\mathbf{x_{nmpc}} = \left[x_D,\,\,V_N,\,\,V_D,\,\,\theta,\,\,q,\,\,\delta_{\text{thrust}},\,\,\delta_e,\,\,\Delta{\delta}_\text{thrust}\,\,\Delta{\delta}_e\right]^\intercal,
\end{equation}

The following optimal control problem is solved at each time step:

\begin{equation}
\begin{split}
\min_{\mathbf{x},\mathbf{u}}\, \frac{H_p}{2}\;\sum_{j = 1}^{j = N+1} \bigg(&\big\Vert \bold{x}_D(j) - \bold{x}_{D_\text{ref}}(j)\big\Vert_{Q_x}^2 + \big\Vert \bold{V}_t(j) - \bold{V}_{t_\text{ref}}(j)\big\Vert_{Q_{VT}}^2 \\
& + \big\Vert \bold{V}_D(j) - \bold{V}_{D_\text{ref}}(j)\big\Vert_{Q_{VD}}^2 + \big\Vert \Delta{\delta}_{\text{thrust}}\big\Vert_{Q_{\text{T}}}^2\\
&  + \big\Vert \Delta{\delta}_e\big\Vert_{Q_{\delta_e}}^2 + \big\Vert q\big\Vert_{Q_{q}}^2 + \big\Vert a_D\big\Vert_{Q_{a}}^2\bigg)\;w(j),
\end{split}
\end{equation}

subject to
\begin{eqnarray}
\left(\frac{t_f-t_0}{2}\right)\mathbf{D}_{j,k}\mathbf{x}_j - \mathbf{\dot{x}}_j &=& 0, \\
\mathbf{x}(j_0) - \mathbf{x}_{\text{dem}}(j_0) &=& 0,\\
\mathbf{x}_{lb}  \leq    \mathbf{x}  \leq  \: \mathbf{x}_{ub},\\
\mathbf{u}_{lb}  \leq    \mathbf{u}  \leq  \: \mathbf{u}_{ub},\\
\Delta\mathbf{u_{\text{lb}}}  \leq    \Delta\mathbf{u}  \leq  \: \Delta\mathbf{u_{\text{ub}}}, 
\label{eq:chap5_6DOF_cons}
\end{eqnarray}

where $V_T$ and $V_{T_\text{ref}}$ are the actual and reference true airspeeds respectively and the state vector $\mathbf{x}$ is defined in \eqref{eqn:xnmpc_thrust}.  $Q_x$, $Q_{VT}$, $Q_{VD}$, $Q_T$, $Q_{\delta_e}$, $Q_q$ and $Q_a$ are diagonal weighting matrices with the following values along the diagonals 10, 5, 5, 0.001, 0.1, 0.01 and 0.01 respectively and the term $w(j)$ are the pseudospectral node weights \cite{RKPaper1}.  The constraints applied are: $x_D$: 1 to 300m, $V_N$: $30$ to $100\,\text{m/s}$, $V_D$: $\pm 3\,\text{m/s}$, $\theta$: None, $q$: None, $\delta_e$: $\pm 20\deg$, $\Delta{\delta}_{\text{thrust}}$: $\pm 6500\,\text{N/s}$ and $\Delta{\delta}_e$: $\pm 200\,\text{deg/s}$.

The lower limit on thrust is $0\,\text{N}$ while the upper limit changes throughout the flight and is set to the maximum value of  thrust based on the height of the aircraft.  Maximum thrust is calculated via equation \eqref{eqn:chap5_thrustModel}.  If a fault has been detected and the fault flag described in section \ref{subsubsec:chap5_FDL} is set to 1 the upper constraint is set to the filter estimate of thrust plus a $2\sigma$ uncertainty.

The following scenarios were designed to test the fault tolerant control system:

\begin{enumerate}[label=\bfseries Scenario \arabic*:, leftmargin = 100pt]
\item no fault case
\vspace{-10pt}
\item engine failure - $65\%$ power loss 30 secs into flight,
\vspace{-10pt}
\item engine failure - $70\%$ power loss 30 secs into flight.
\end{enumerate}

Note: all test runs take into account the effect of wind.

\begin{figure}[H]
\hspace{-0.5in}
\includegraphics[scale=0.4]{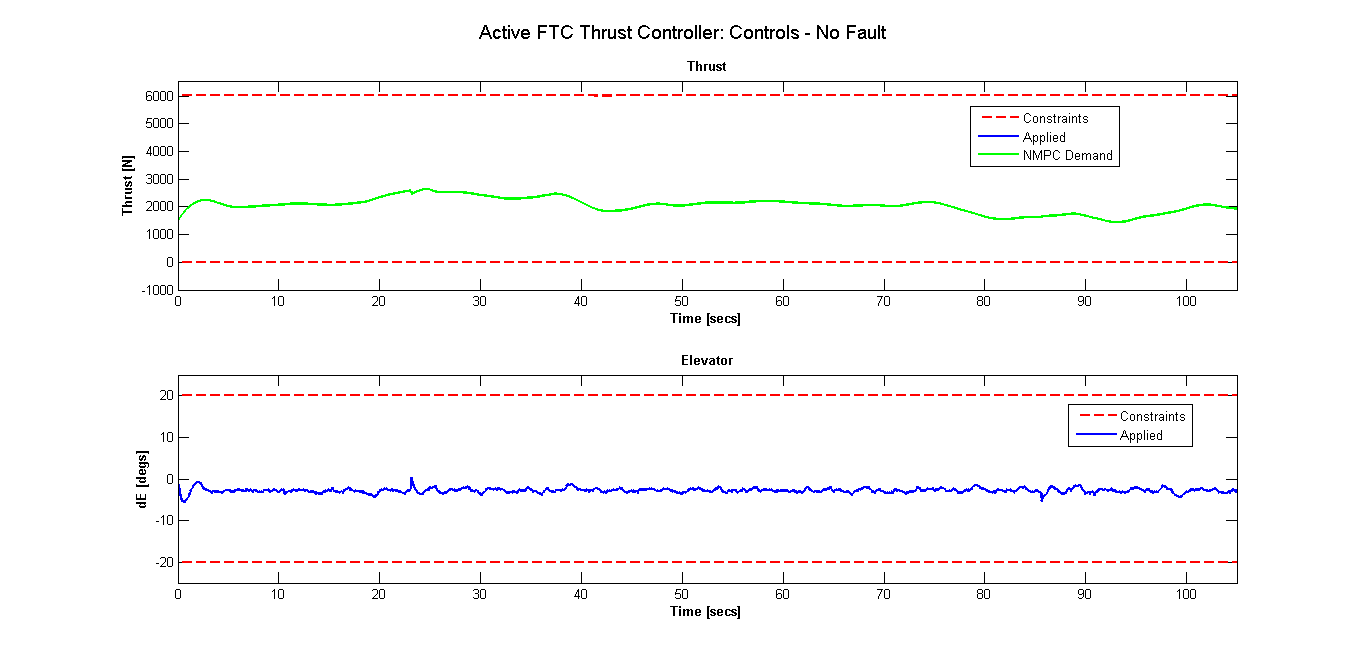} 
%\vspace{-0.5in}
\caption{Active FTC Thrust Controller: Control Inputs - Scenario 1: No Fault Case}
\label{fig:chap5_activeFTC_TC_Controls_NF}
\end{figure}

\begin{figure}[H]
\hspace{-0.5in}
\includegraphics[scale=0.4]{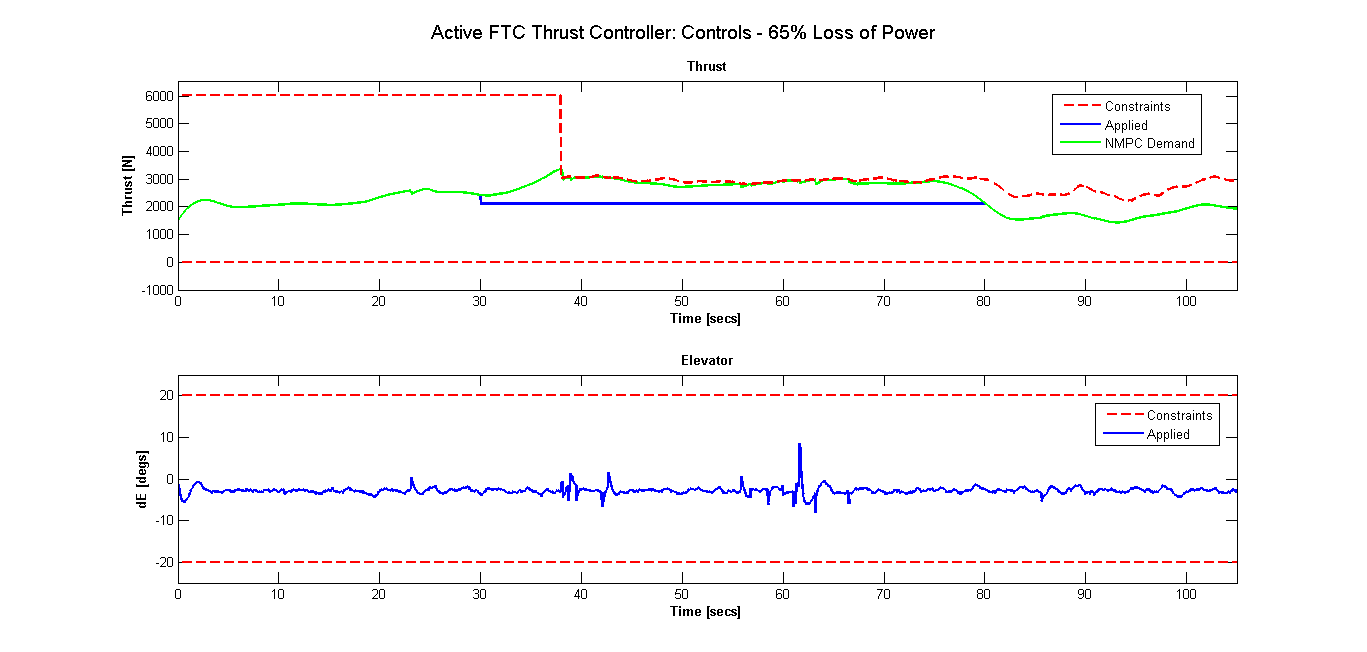} 
%\vspace{-0.5in}
\caption{Active FTC Thrust Controller: Control Inputs - Scenario 2: $65\%$ Loss of Power Case}
\label{fig:chap5_activeFTC_TC_Controls_35T}
\end{figure}

\begin{figure}[H]
\hspace{-0.5in}
\includegraphics[scale=0.4]{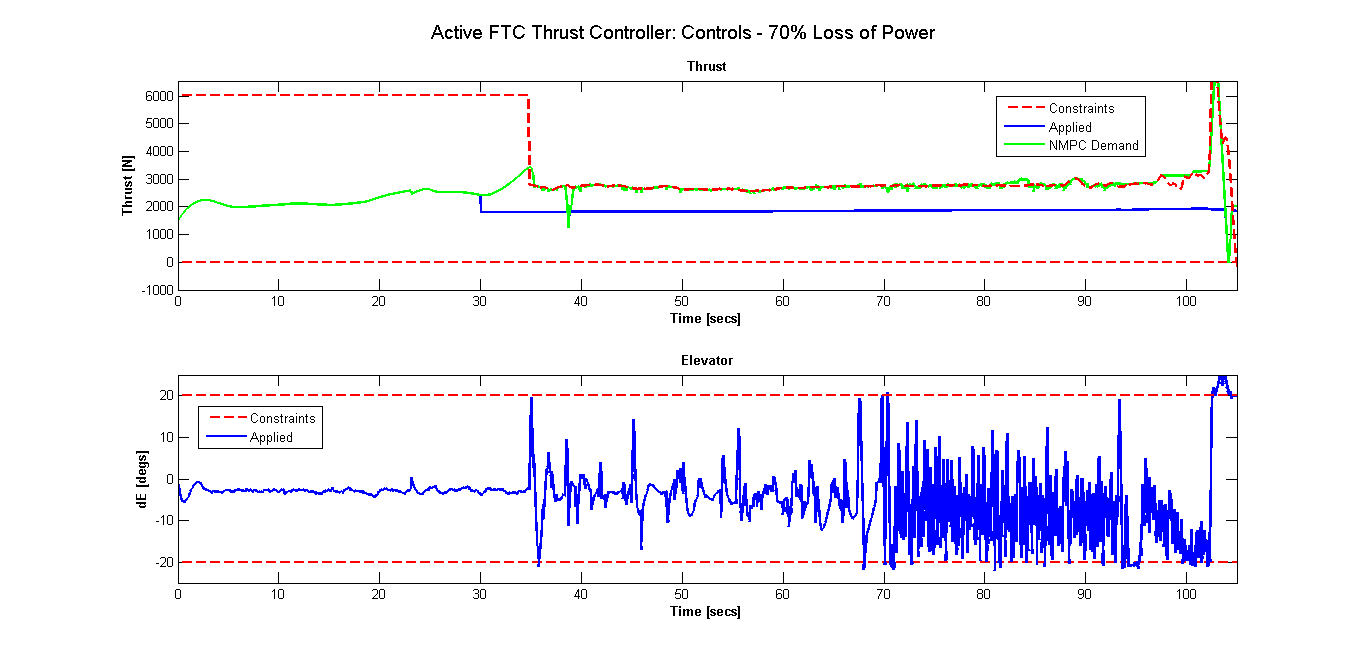} 
%\vspace{-0.5in}
\caption{Active FTC Thrust Controller: Control Inputs - Scenario 3: $70\%$ Loss of Power Case}
\label{fig:chap5_activeFTC_TC_Controls_30T}
\end{figure}

Figures \ref{fig:chap5_activeFTC_TC_Controls_NF}, \ref{fig:chap5_activeFTC_TC_Controls_35T} and \ref{fig:chap5_activeFTC_TC_Controls_30T} show the control inputs for scenarios 1, 2 and 3 respectively.  The results for the thrust show a dip in the constraint value for the upper thrust limits for scenarios 2 and 3 soon after 30 secs.  This indicates that the fault was correctly identified and the NMPC was reconfigured with the information provided by the FDI filter.  The uncertainty bounds in both figures are slightly higher than the actual thrust applied due to the addition of the $2\sigma$ uncertainty.  Other values of $\sigma$ were found to cause the controller and hence the filter to diverge.  Although the estimate is slightly above the actual it is still in the vicinity of the actual thrust level and prompts the controller to allocate more control authority to the other available actuators.  The results show that compared to the no fault case once a fault occurs the elevator activity increases as the power decreases.  Also the more severe the fault the faster the detection time.  This is evident from the fact that the fault is detected earlier in the $70\%$ power loss case compared to the $65\%$ loss of power case.  

\begin{figure}[H]
\hspace{-0.5in}
\includegraphics[scale=0.4]{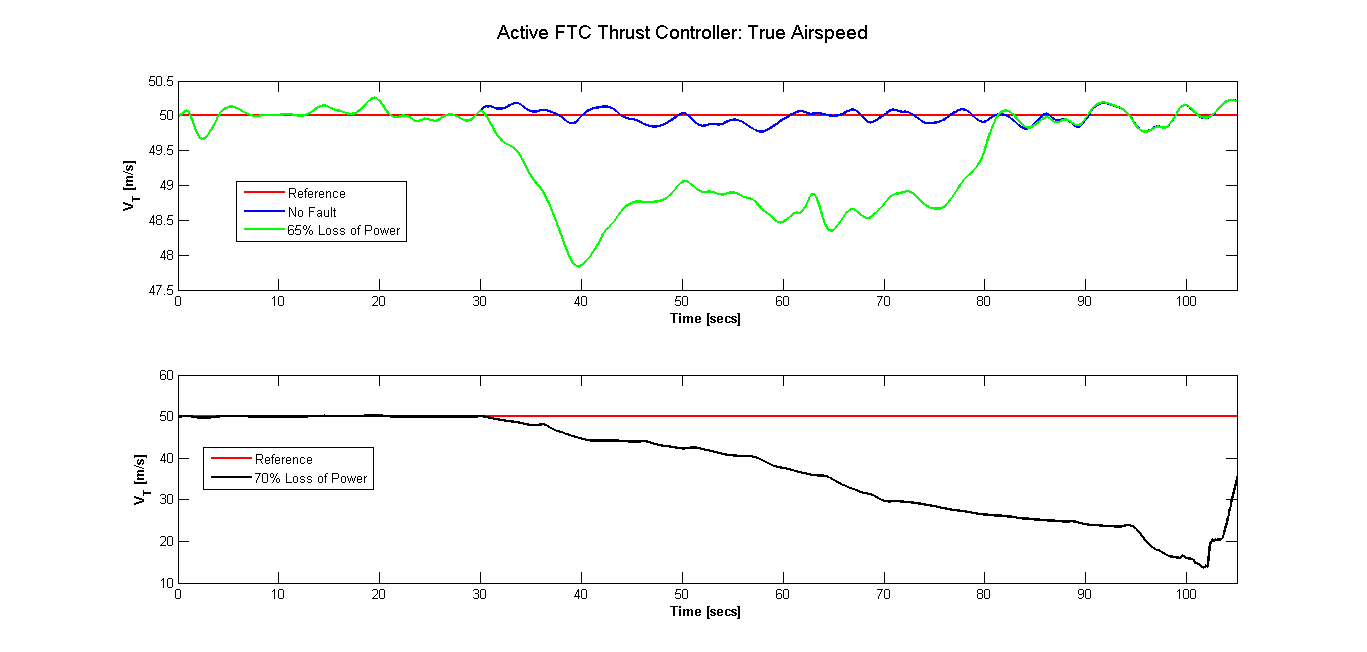} 
%\vspace{-0.5in}
\caption{Active FTC Thrust Controller: True Airspeeds}
\label{fig:chap5_activeFTC_TC_Vt}
\end{figure}

Figure \ref{fig:chap5_activeFTC_TC_Vt} shows the true airspeed of the air vehicle.  A true airspeed of 50m/s was demanded by the aircraft.  In the case where there is $65\%$ loss of power the aircraft is unable to maintain the demanded true airspeed during straight and level flight.  Once the aircraft begins to descend the demanded true airspeed is recovered.  However in the case of $70\%$ power loss there is not enough power to maintain the demanded airspeed.  Once the fault occurs the airspeed begins to drop and reaches stall causing the aircraft to lose control.
  
\begin{figure}[H]
\hspace{-0.5in}
\includegraphics[scale=0.4]{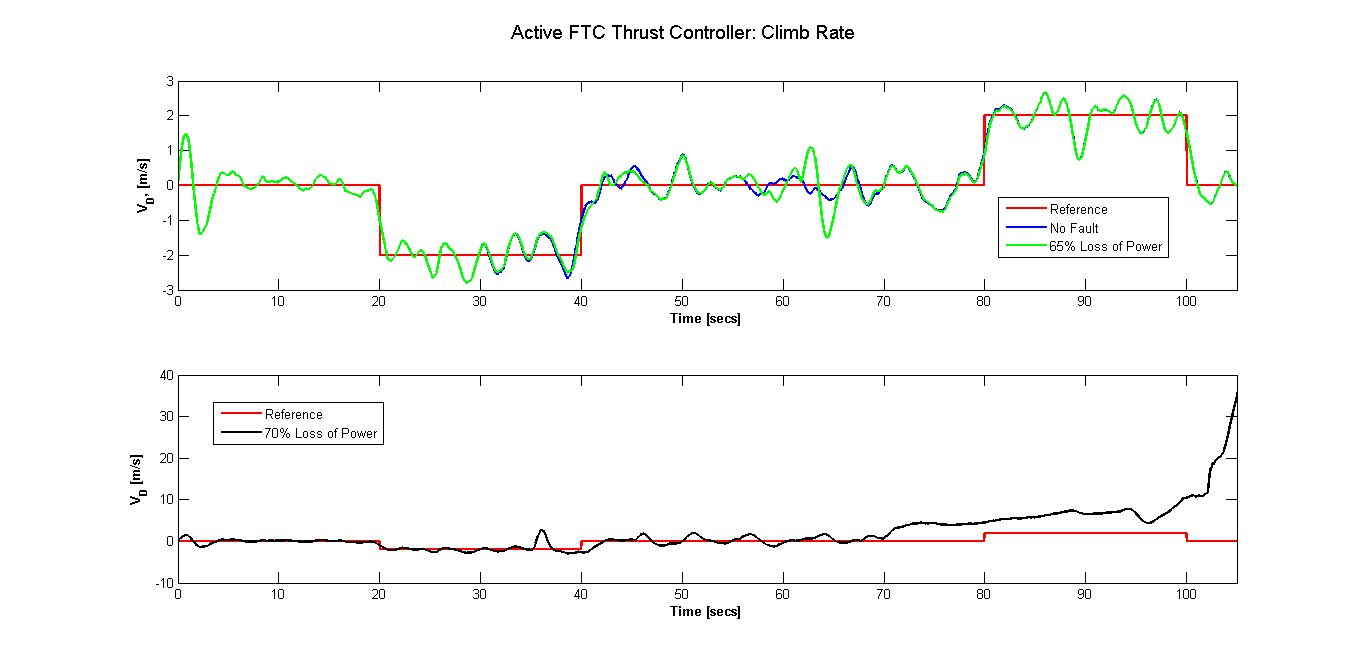} 
%\vspace{-0.5in}
\caption{Active FTC Thrust Controller: Climb Rates}
\label{fig:chap5_activeFTC_TC_Vd}
\end{figure}

The vertical speed (also known as climb rate) response is given in figure \ref{fig:chap5_activeFTC_TC_Vd}.  In the $65\%$ power loss case the response is very similar to the no fault response.  A $70\%$ loss in power results in the aircraft being unable to maintain speed and it descends to the ground.

\begin{figure}[H]
\hspace{-0.5in}
\includegraphics[scale=0.4]{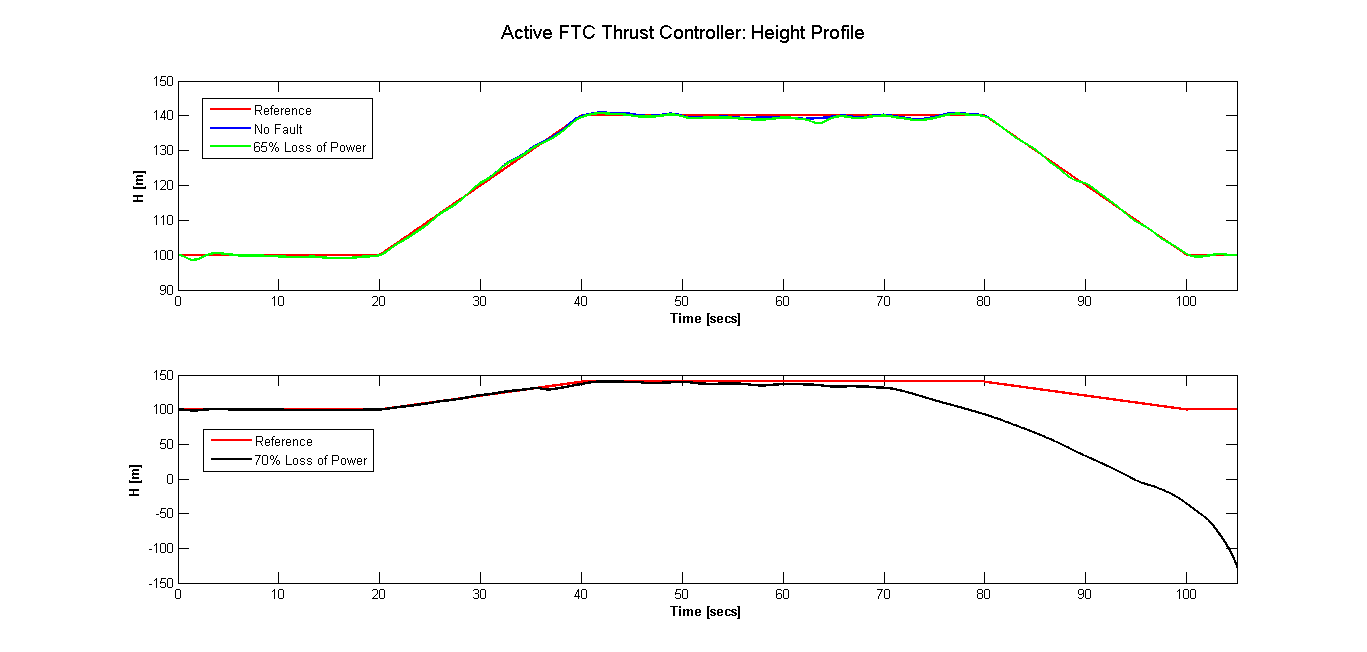} 
%\vspace{-0.5in}
\caption{Active FTC Thrust Controller: Height Profiles}
\label{fig:chap5_activeFTC_TC_Height}
\end{figure}

The height profiles given in figure \ref{fig:chap5_activeFTC_TC_Height} show that even with a $65\%$ loss in power the aircraft is capable of maintaining the reference trajectory.  However when the power decreases by another $5\%$ the aircraft completes the climb to the highest demanded altitude but begins to descend half way through straight and level flight.

In summary, we can see that the results obtained successfully demonstrate the application of the active fault tolerant flight control system design.  The control system is able to detect an engine fault within 2-3 seconds of the fault occurring, enabling reconfiguration of the NMPC controller to allow reallocation of control authority to maintain the aircraft on the demanded flight path within the aircraft limits.  The controller works hard to achieve the demands, however in the event where this is impossible this information can be used to bring the aircraft back safely to the ground.

\section{Conclusion}\label{section:conclusion}
In this paper we effectively demonstrated the application of our active FTC system design to flight control.  The FTC system comprised of an NMPC controller integrated with a UKF filter for fault detection.  This is the first time such a system has been applied within the context of fault tolerant flight control.  To assist the research a generic aircraft model was developed and the active FTC system was applied to the 6DoF aircraft model.  Research into the application of the FTC system on a full 6Dof aircraft model identified a number of areas for further research particularly in the design of the FDI system.  The FTC system was then applied to the longitudinal motion of the aircraft where engine failure scenarios were simulated.  The results obtained show that the FTC system successfully identified the fault within seconds of occurrence and re-allocated the control authority to healthy actuators based upon up to date fault information.  

\section*{Appendix}
\subsubsection{Force and Moment Coefficients}\label{subsubsec:chap5_FandMCoeffs}
For the F-4 at $\alpha\leq 15^0$ the non-dimensional force and moment coefficients are given by:
\begin{equation}
\begin{split}
CX ={}& -0.0434 + 2.93 \times 10^{-3}\alpha + 2.53\times 10^{-5}\beta^2-1.07\times 10^{-6}\alpha \beta^2 + 9.5 \times 10^{-4}\delta_e\\
& -8.5\times 10^{-7}\delta_e \beta^2 + \left(\frac{180q\bar{c}}{\pi 2 V_t}\right)\left(8.73\times 10^{-3} + 0.001\alpha - 1.75 \times 10^{-4} \alpha^2 \right),
\end{split} \label{eqn:chap5_CX}
\end{equation}

\begin{equation}
\begin{split}
CY ={}& -0.012\beta + 1.55\times 10^{-3}\delta_r - 8\times 10^{-6}\delta_r\alpha\\
& +\left(\frac{180b}{\pi 2 V_t}\right)\left(2.25 \times 10^{-3} p + 0.0117r-3.67\times 10^{-4}r\alpha + 1.75\times 10^{-4}r\delta_e\right),
\end{split}\label{eqn:chap5_CY}
\end{equation}

\begin{equation}
\begin{split}
CZ ={}& -0.131 - 0.0538\alpha - 4.76\times 10^{-3}\delta_e - 3.3\times 10^{-5}\delta_e\alpha-7.5\times 10^{-5}{\delta_a}^2\\
& +\left(\frac{180q\bar{c}}{\pi 2 V_t}\right)\left(-0.111+5.17 \times 10^{-3} \alpha - 1.1\times 10^{-3}{\alpha}^2\right),
\end{split}\label{eqn:chap5_CZ}
\end{equation}

\begin{equation}
\begin{split}
C_l = {}& -5.98\times 10^{-4}\beta - 2.83 \times 10^{-4} \alpha \beta + 1.51\times 10^{-5}\alpha^2 \beta\\
& -\delta_a \left(6.1\times 10^{-4} + 2.5\times 10^{-5}\alpha - 2.6\times 10^{-6}\alpha^2\right)\\
&+\delta_r\left(-2.3\times 10^{-4}+4.5\times 10^{-6}\alpha\right)\\
& +\left(\frac{180b}{\pi 2 V_t}\right)\left(-4.2\times 10^{-3}p-5.24\times 10^{-4}p\alpha + 4.36\times 10^{-5}p\alpha^2 \right. \\
& \left. + 4.36\times 10^{-4}r + 1.05\times 10^{-4}r\alpha + 5.24\times 10^{-5}r\delta_e\right),
\end{split}\label{eqn:chap5_Cl}
\end{equation}

\begin{equation}
\begin{split}
C_m = {}& -6.61\times 10^{-3} - 2.67 \times 10^{-3} \alpha -6.48\times 10^{-5}\beta^2\\
& -2.65\times 10^{-6}\alpha\beta^2 - 6.54\times 10^{-3}\delta_e - 8.49\times 10^{-5}\delta_e\alpha\\
& + 3.74\times 10^{-6}\delta_e\beta^2 - 3.5\times 10^{-5}{\delta_a}^2\\
&+ \left(\frac{180q\bar{c}}{\pi 2 V_t}\right)\left(-0.0473-1.57\times 10^{-3}\alpha\right)+\left(x_{c.g.ref}-x_{c.g}\right)C_Z,
\end{split}\label{eqn:chap5_Cm}
\end{equation}

\begin{equation}
\begin{split}
C_n = {}& 2.28\times 10^{-3}\beta + 1.79\times 10^{-6}\beta^3 + 1.4\times 10^{-5}\delta_a\\
& + 7.0\times 10^{-6}\delta_a\alpha - 9.0\times 10^{-4}\delta_r + 4.0 \times 10^{-6}\delta_r\alpha\\
& + \left(\frac{180b}{\pi 2 V_t}\right)\left(-6.63\times 10^{-5}p-1.92\times 10^{-5}p\alpha + 5.06\times 10^{-6}p\alpha^2 \right.\\
&\quad\quad\quad\quad\quad\quad \left. -6.06\times 10^{-3}r-8.73\times 10^{-5}r\delta_e+8.7\times 10^{-6}r\delta_e\alpha\right)\\
&- \left(\frac{\bar{c}}{b}\right)\left(x_{c.g.ref}-x_{c.g}\right)C_Z.
\end{split}\label{eqn:chap5_Cn}
\end{equation}

\subsubsection{Thrust Model} 
The following thrust model \cite{bryson1999dynamic} is used in this work:

\begin{align} \label{eqn:chap5_thrustModel}
\begin{split}
h_T = {}&\frac{H}{3048},
\end{split}\\
\nonumber \\
\begin{split}
T_{max} ={}&((30.21-0.668\,{h_T}-6.877\,{h_T}^2+1.951\,{h_T}^3-0.1512\,{h_T}^4)\\
    &+ \left(\frac{Vt}{v_s}\right)(-33.8+3.347\,{h_T}+18.13\,{h_T}^2-5.865\,{h_T}^3+0.4757\,{h_T}^4)\\
     &+\left(\frac{Vt}{v_s}\right)^2(100.8-77.56\,{h_T}+5.441\,{h_T}^2+2.864\,{h_T}^3-0.3355\,{h_T}^4)\\
      &+\left(\frac{Vt}{v_s}\right)^3(-78.99+101.4\,{h_T}-30.28\,{h_T}^2+3.236\,{h_T}^3-0.1089\,{h_T}^4)\\ 
    &+\left(\frac{Vt}{v_s}\right)^4(18.74-31.6\,{h_T}+12.04\,{h_T}^2-1.785\,{h_T}^3+0.09417\,{h_T}^4))\frac{4448.22}{20},
\end{split}\\
\nonumber\\
T ={}& T_{max}\,\delta_{th}, 
\end{align}

where $v_s$ is the speed of sound, $340.3\,m/s$, $H$ the height or the $-x_D$ position of the aircraft, and $\delta_{th}$.

\end{document}